\documentclass[article]{amsart}

\usepackage{graphics}
\usepackage[all]{xy}

\usepackage{amssymb,amsmath}
\usepackage{psfrag}

\usepackage{mathrsfs} 
\usepackage{amsfonts}
\usepackage{amscd}

\usepackage{graphicx}
\usepackage{epstopdf}
\usepackage{tikz}
\usetikzlibrary{arrows,automata}

\usepackage{graphicx,epsfig,float}

\usepackage{url}

\makeatletter
\def\url@leostyle{%
  \@ifundefined{selectfont}{\def\UrlFont{\sf}}{\def\UrlFont{\small\ttfamily}}}
\makeatother
\newtheorem{thm}{Theorem}[section]
\newtheorem{lem}[thm]{Lemma}

\newtheorem{cor}[thm]{Corollary}
\newtheorem{fact}[thm]{Fact}

\newtheorem{clm}[thm]{Claim}

\newtheorem{defn}[thm]{Definition}

\newtheorem{nrmk}[thm]{Remark}

\newtheorem{propt}[thm]{Properties}

\newcommand{\pf}{{\bf Proof. }}

\renewcommand{\tilde}{\widetilde}
\renewcommand{\bar}{\overline}

\newcommand{\bJ}{{\mathbf J}}

\newcommand{\maI}{\mathcal I}
\newcommand{\maJ}{\mathcal J}
\newcommand{\maA}{\mathcal A}

\newcommand{\maD}{\mathcal D}
\newcommand{\maK}{\mathcal K}
\newcommand{\maY}{\mathcal Y}
\newcommand{\maL}{\mathcal L}
\newcommand{\maR}{\mathcal R}
\newcommand{\maU}{\mathcal U}
\newcommand{\maX}{\mathcal X}

\newcommand{\od}{^{\mathrm{op}}}

\newcommand{\F}{\mbox{$\mathcal{F}$}}

\newcommand{\psh}{\mathrm{Psh}}
\newcommand{\sh}{\mathrm{Sh}}

\newcommand{\ob}{\mathrm{Ob}}
\newcommand{\op}{\mathrm{Op}}

\newcommand{\Ho}{\mathrm{Hom}}

\def\de{\mathrm{def}}

\def\Aut{\mathrm{Aut}}

\newcommand{\lcsh}{\mathrm{LCSh}}
\newcommand{\csh}{\mathrm{CSh}}

\def\fct{\mathrm{Fct}}

\begin{document}

\title[]{Fundamental group in o-minimal structures with definable Skolem functions}

\author{Bruno Dinis}

\address{ Departamento de Matem\'atica\\
Faculdade de Ci\^encias,  Universidade de Lisboa\\
Campo Grande, Edif\'icio C6\\ 
P-1749-016 Lisboa, Portugal}

\email{bmdinis@fc.ul.pt}

\author {M\'ario J. Edmundo}

\address{ Departamento de Matem\'atica\\
Faculdade de Ci\^encias,  Universidade de Lisboa\\
Campo Grande, Edif\'icio C6\\ 
P-1749-016 Lisboa, Portugal}

\email{mjedmundo@fc.ul.pt}

\author{Marcello Mamino}

\address{Universit\`a di Pisa \\
Largo Bruno Pontecorvo 5 \\
56127 Pisa, Italia}

\email{marcello.mamino@dm.unipi.it}

\date{\today}
\thanks{The first author would like to acknowledge the support
of Centro de Matem\'atica e Aplica\c{c}\~oes Fundamentais com Investiga\c{c}\~ao Operacional [project
UID/MAT/04561/2019].\newline
 }


\begin{abstract}
In this paper we work in an arbitrary o-minimal structure with definable Skolem functions and we prove that definably connected, locally definable manifolds are uniformly definably path connected, have an admissible cover by definably simply connected, open definable subsets and, definable paths and definable homotopies on such locally definable manifolds can be lifted to locally definable covering maps. These properties  allows us to obtain the main properties of the general o-minimal fundamental group, including:  invariance and comparison results;  existence of  universal locally definable covering maps;  monodromy equivalence for locally constant o-minimal sheaves - from which one obtains, as in algebraic topology, classification results for locally definable covering maps, o-minimal Hurewicz and Seifert - van Kampen theorems. 
\end{abstract}

\maketitle

\begin{section}{Introduction}\label{section intro}
In this paper we work in an  arbitrary o-minimal structure ${\mathbb M}=(M,<, (c)_{c\in {\mathcal C}}, $ $ (f)_{f\in  {\mathcal F}}, (R)_{R\in {\mathcal R}})$ which we assume to have definable Skolem functions.  We are interested in developing algebraic topology tools for objects definable  in ${\mathbb M},$ more specifically a suitable fundamental group functor. One should point out here that, although the objects definable in  ${\mathbb M}$ have a topology induced by the order in $M,$  if ${\mathbb M}$ is non-archimedean, then all  such objects (except the ones of dimension zero) are totally disconnected topological spaces and so the topological fundamental group is of no use.

In the case ${\mathbb M}$ expands a real closed field $(M,<,0,1, +, \cdot),$ so this includes the semi-algebraic case, a suitable o-minimal fundamental group functor is already known and was studied in, for example, \cite{dk4}, \cite{dk5}, \cite{bo1}, \cite{bo2}, \cite{eo}, \cite{bao1} and \cite{bao2}, where its main properties were proved and applications to the theory of definable groups were obtained. These main properties include: (i) finite generation; (ii) invariance when going to o-minimal expansions of ${\mathbb M}$ or to bigger models of the first order theory of ${\mathbb M}$ and (iii) coincidence with the topological fundamental group when $M={\mathbb R}.$ Often, in this case, the proofs of these properties relied on the o-minimal triangulation theorem (\cite{vdd}), a generalization of the semi-algebraic and of the semi-analytic triangulation theorems (\cite{dk2}, \cite{Gi}, \cite{Lo} and \cite{Hi}). 

In the case ${\mathbb M}$ expands an ordered group $(M,<,0, + )$ a suitable o-minimal fundamental group functor was studied in the papers \cite{el}, \cite{es}, \cite{EdPa} and \cite{eep}, where besides the properties mentioned above, the following were also proved: (iv) the  existence of  universal locally definable covering maps; (v)  monodromy equivalence for locally constant o-minimal sheaves - from which one obtains, as in algebraic topology: (vi) classification results for locally definable covering maps; (vii) o-minimal Hurewicz and Seifert - van Kampen theorems. 

In the case ${\mathbb M}$ is an arbitrary o-minimal structure with definable Skolem functions, following M. Mamino's ideas, a general o-minimal fundamental group was introduced in the paper \cite{Edal17}, where Pillay's conjecture for definably compact groups was obtained in the general case. However, in that paper, the properties mentioned above were only proved for the o-minimal $\bJ$-fundamental group which is the relativization of the general o-minimal fundamental group to a cartesian product  $\bJ=\Pi _{i=1}^mJ_i$  of   definable group-intervals $J_i=\langle (-_ib_i, b_i), $ $0_i, +_i, -_i,  <\rangle $.  The goal of this paper is to show, as conjectured in \cite{Edal17}, all of the above properties for the general o-minimal fundamental group in the more general setting of o-minimal structures with definable Skolem functions.

The main technical result of the paper, obtained in Subsection \ref{subsec main prop}, says that definably connected, locally definable manifolds are uniformly definably path connected, have an admissible cover by definably simply connected, open definable subsets and, definable paths and definable homotopies on such locally definable manifolds can be lifted to locally definable covering maps. See Properties \ref{propt main new paths}.   The results (i), (ii), (iii) and (iv) mentioned above are obtained, in the general case of ${\mathbb M}$  an arbitrary o-minimal structure with definable Skolem functions,  in Subsection \ref{subsec inv main}, and as explained in Subsection \ref{subsec main prop}, they follow from Properties \ref{propt main new paths} in exactly the same way as  in \cite{eep} for o-minimal expansions of ordered groups.   In Subsection \ref{subsec inv main} we also show that if  $X$  is a definably connected locally definable manifold with definable charts in $\bJ,$ then the o-minimal  $\bJ$-fundamental group of $X$  is isomorphic to the o-minimal fundamental group of $X$ in ${\mathbb M}.$ Result (v) and its consequences (vi) and (vii) are mentioned in Subsection \ref{subsection monodromy} after some background is recalled.

It might be possible that part of this theory may be developed in a general o-minimal structure without the definable Skolem functions assumption. However, in the applications that we have in mind, all the structures have definable Skolem functions. 
It is even possible that under the assumption of the existence of definable Skolem functions some arguments may be simplified. Indeed, due to the trichotomy theorem it might be possible that there is a definable family of groups covering all the points of the domain of the structure. In that case the arguments in Subsection \ref{subsec main prop} would be closer to those of the o-minimal $\bJ$ fundamental group treated in \cite{Edal17}. We leave these issues for future work.

Regarding the applications, consider an algebraically closed valued field with a nontrivial valuation and consider the definable sets induced in the value group with a point at infinity added (the valuation of $0$), $\Gamma_{\infty}$. Adding a copy of the value group $\Gamma$ to $\Gamma_{\infty}$ we obtain an o-minimal structure $\Sigma$ with definable Skolem functions. Understanding the definable topology/definable algebraic topology of definable subsets of $(\Gamma_{\infty})^n$ can have applications in non-Archimedean tame topology by the main theorem of Hrushovski and Loeser (\cite[Theorem~11.1.1]{HL16}). For instance, Theorem~\ref{thm continuous-choice} below is used in an essential way in developing the cohomology in the non-Archimedean tame setting \cite{EKY(ta)}. Moreover, we expect that the main results about the fundamental group developed here or adaptations of those might give further applications.\\
\end{section}

\begin{section}{Preliminaries}\label{section prelimp}
In this section we recall the notion of locally definable manifolds and locally definable covering maps and the general o-minimal fundamental group from \cite{Edal17}.

Before we start, recall that an o-minimal structure ${\mathbb M}$ has {\it definable Skolem functions}  if and only if for every uniformly definable family $\{X_t\}_{t\in T}$ of nonempty definable subsets of some $M^k,$ there is a definable function $h:T\to M^k$ such that:
\begin{itemize}
\item[-]
$h(t)\in X_t$ for all $t\in T.$
\end{itemize}

\begin{subsection}{Locally definable manifolds  and covering maps}\label{subsection ldef man covers}
Here we recall the definition of the category of locally definable manifolds with continuous locally definable maps and the  notion of locally definable covering maps.\\

A {\it locally definable manifold (of dimension $n$)} is a triple $(S,(U_i,\theta _i)_{i\leq \kappa})$ 
where:

 \begin{itemize}
\item[$\bullet  $]
$S=\bigcup _{i\leq \kappa}U_i$;

\item[$\bullet  $]
each $\theta  _i:U_i\rightarrow  M^{n}$ is an injection such that $\theta _i(U_i)$ is an open definable subset of $M^{n}$;

\item[$\bullet  $]
for all $i, j$, $\theta  _i(U_i\cap U_j)$ is an open definable subset of  $\theta _i(U_i)$ and  the transition maps $\theta _{ij}:\theta  _i(U_i\cap U_j)\rightarrow  \theta  _j(U_i\cap U_j):x\mapsto \theta  _j(\theta  _i^{-1}(x))$ are definable homeomorphisms.
\end{itemize}
We call the $(U_i, \theta _i)$'s the {\it definable charts of $S$}. If $\kappa <\aleph _0 $ then $S$ is a {\it definable manifold}. 

A locally definable manifold $S$ is equipped with the topology such that  a subset $U$ of $S$ is  open  if and only if for each $i$, $\theta  _i(U\cap U_i)$ is an open definable subset of $\theta  _i(U_i)$.

We  say that a subset $A$ of $S$ is \emph{definable} if and only if there is a finite $I_0\subseteq \kappa $ such that $A\subseteq \bigcup _{i\in I_0}U_i$ and for each $i\in I_0$,  $\theta _i(A\cap U_i)$ is a definable subset of $\theta  _i(U_i)$. A subset $B$ of $S$ is \emph{locally definable} if and only if for each $i$,  $B\cap U_i$ is a definable subset of $S$.   We say that  a locally definable manifold $S$ is {\it definably connected} if it is not the disjoint union of two open and closed locally definable subsets. 

If ${\mathcal U}=\{U_{\alpha }\}_{\alpha \in I}$ is a cover of $S$ by open locally definable subsets, we say that ${\mathcal U}$  is  {\it admissible} if for each $i\leq \kappa  $, the cover $\{U_{\alpha }\cap U_i\}_{\alpha \in I}$ of $U_i$ admits a finite subcover. If ${\mathcal V}=\{V_{\beta }\}_{\beta \in J}$ is another cover of $S$ by open locally definable subsets, we say that ${\mathcal V}$ {\it refines} ${\mathcal U}$, denoted by ${\mathcal V}\leq {\mathcal U},$ if there is a map $\epsilon :J\rightarrow  I$ such that $V_{\beta }\subseteq U_{\epsilon (\beta )}$ for all $\beta \in J$.

A map $f:X\to  Y$ between locally definable manifolds with definable charts $(U_i,\theta _i)_{i\leq \kappa _X}$ and $(V_j,\delta _j)_{j\leq \kappa _Y}$ respectively is {\it a locally definable map} if for every finite $I\subseteq \kappa _X$ there is a finite $J\subseteq \kappa _Y$ such that:
\begin{itemize}
\item[$\bullet$]  $f(\bigcup _{i\in I}U_i)\subseteq \bigcup _{j\in J} V_j$;

\item[$\bullet$] the  restriction $f_{|}:\bigcup _{i\in I}U_i \to  \bigcup _{j\in J} V_j$ is a definable map between definable manifolds, i.e., for each $i\in I$ and every $j\in J$,
    $\delta _j\circ f\circ \theta _i^{-1}:\theta _i(U_i)\to  \delta _j(V_j)$ is a definable map between definable sets.
\end{itemize}
Thus we have the category of locally definable manifolds with locally definable continuous maps.\\

\begin{defn}
{\em
Given a  definably connected locally definable manifold $S$, a locally definable manifold $X$ and an admissible cover ${\mathcal U}=\{U_{\alpha }\}_{\alpha \in I}$ of $S$ by open definable subsets, we say that a continuous surjective locally definable map $p_X:X\rightarrow S$ is a {\it locally definable covering map trivial over ${\mathcal U}=\{U_{\alpha }\}_{\alpha \in I}$} if the following hold:
\begin{enumerate}
\item[$\bullet $]
$p_X^{-1} (U_{\alpha })=\bigsqcup _{i\leq \lambda }U_{\alpha }^i$ a disjoint union of open definable subsets of $X;$

\item[$\bullet  $]
each $p_{X|U_{\alpha }^i}:U_{\alpha } ^i\rightarrow U_{\alpha }$ is a definable homeomorphism.
\end{enumerate}
A {\it locally definable covering map} $p_X:X\rightarrow S$ is a locally  definable covering map trivial over some admissible cover ${\mathcal U}=\{U_{\alpha }\}_{\alpha \in I}$ of $S$ by open definable subsets.

We say that two locally definable covering maps $p_X:X\rightarrow S$ and $p_Y:Y\rightarrow S$ 
are {\it locally definably homeomorphic} if there is a locally definable homeomorphism $F:X\rightarrow Y$ such that:

\begin{enumerate}
\item[$\bullet $]
$p_X=p_Y\circ F.$
\end{enumerate}

A locally definable covering map $p_X:X\to  S$ is {\it trivial} if it is locally definably homeomorphic to a locally definable covering map $S\times M\to S :(s,m)\mapsto s$ for some  set $M.$\\
}
\end{defn}

Let $p_Y:Y\to  T$ be a locally definable covering map, $X$ be a  locally definable manifold  and let $f:X\to  T$ be a locally definable  map.  A {\it lifting of $f$} is a continuous map $\tilde{f}:X\to  Y$  such that $p_Y\circ \tilde{f}=f$. Note that a lifting of a continuous  locally definable  map need not be a 
locally definable  map. However, if $X$ is definably connected, then any two continuous locally definable liftings which coincide in a point must be equal \cite[Lemma 2.8]{eep}.\\

\end{subsection}

\begin{subsection}{A general o-minimal fundamental group}\label{subsection new fund grp}
Here we recall the definition and the basic properties of the  o-minimal fundamental group  in arbitrary o-minimal structures from \cite{Edal17}. \\

\begin{defn}
{\em
By a {\it basic $d$-interval}, short for {\it basic directed interval},  we mean a tuple 
$$\maI=\langle [a,b], \langle 0_{\maI}, 1_{\maI}\rangle \rangle$$
where $a,b \in M$ with $a<b$ and $\langle  0_{\maI}, 1_{\maI}\rangle \in \{\langle a, b\rangle , \langle b, a\rangle \}.$ The {\it domain} of $\maI $ is $[a,b]$ and the {\it direction} of $\maI $ is $\langle 0_{\maI}, 1_{\maI}\rangle.$ The {\it opposite} of  $\maI $ is the basic $d$-interval
$$\maI \od =\langle [a,b], \langle 0_{\maI \od}, 1_{\maI \od} \rangle \rangle $$
with the same domain and opposite direction $\langle 0_{\maI \od}, 1_{\maI \od}\rangle =\langle 1_{\maI }, 0_{\maI} \rangle .$

If $\maI _i=\langle [a_i, b_i], \langle 0_{\maI _i}, 1_{\maI _i}\rangle \rangle$ are basic $d$-intervals, for $i=1, \ldots , n,$ we define the {\it $d$-interval}, short for {\it directed interval}, $ \maI _1\wedge \dots \wedge \maI _n,$  whose {\it domain } is the set 
\[
[a_1, b_1]\wedge \dots \wedge [a_n, b_n]:= \raisebox{1ex}{$\mathsurround=0pt\displaystyle \bigsqcup_ i\{c_i\}\times [a_i, b_i]$}
\Big/ \raisebox{-1ex}{$\mathsurround=0pt\displaystyle\sim$}
\]  
where $c_1,\dots , c_n$ are $n$ distinct points of~$M$ and $\sim $ is the equivalence relation defined by  $(c_i,1_{\maI _i})\sim(c_{i+1}, 0_{\maI _{i+1}})$ for each~$i=1,\dots , n-1$ and identity elsewhere.   The  {\it direction} of $\maI _1\wedge \dots \wedge \maI _n $ is $\langle 0_{\maI _1\wedge \dots \wedge \maI _n}, 1_{\maI _1\wedge \dots \wedge \maI _n} \rangle $ where 
 $0_{\maI _1\wedge \dots \wedge \maI _n}=\langle c_1, 0_{\maI _1} \rangle $ and $1_{\maI _1\wedge \dots \wedge \maI _n}=\langle c_n, 1_{\maI _n}\rangle .$

The {\it opposite} of  $\maI _1\wedge \dots \wedge \maI _n $ is the $d$-interval
$$(\maI _1\wedge \dots \wedge \maI _n ) \od =\langle [a_1, b_1]\wedge \dots \wedge [a_n, b_n], \langle 0_{(\maI _1\wedge \dots \wedge \maI _n) \od}, 1_{(\maI _1\wedge \dots \wedge \maI _n) \od} \rangle \rangle $$
with the same domain and opposite direction 
\[
\langle 0_{(\maI _1\wedge \dots \wedge \maI _n) \od}, 1_{(\maI _1\wedge \dots \wedge \maI _n) \od}\rangle =\langle 1_{\maI _1\wedge \dots \wedge \maI _n }, 0_{ \maI _1\wedge \dots \wedge \maI _n} \rangle .
\]
}
\end{defn}

\begin{fact}\label{fact op}
{\em
If $\maI _i=\langle [a_i, b_i], \langle 0_{\maI _i}, 1_{\maI _i}\rangle \rangle$ are basic $d$-intervals, for $i=1, \ldots , n,$ then $(\maI _1\wedge \dots \wedge \maI _n ) \od = \maI _n \od\wedge \dots \wedge \maI _1 \od .$
}
\end{fact}

Below, for the notion of {\it definable space} we refer the reader to \cite[page 156]{vdd}. We have:

\begin{fact}\cite[Lemma 2.5]{Edal17}\label{fact d-int space}
Let $\maI =\langle I, \langle 0_{\maI}, 1_{\maI}\rangle \rangle $ be a $d$-interval. Then the domain $I$ of $\maI$ is a Hausdorff, definably compact, definable space  of dimension one which is equipped with a definable total order $<_{\maI}.$
\end{fact}




Due to Fact \ref{fact d-int space}, below we will identify a $d$-interval $\maI =\langle I, \langle 0_{\maI}, 1_{\maI}\rangle \rangle $ with its domain equipped with the definable total order $<_{\maI}.$ In particular, since the domain $I$ of $\maI \od$ is a definable space of dimension one which is equipped with the definable total order $>_{\maI},$ we  have an
order reversing definable homeomorphism (with respect to the topologies given by the orders)
$$o_{\maI}: \maI \to \maI \od$$
given by the identity on the domain. \\

Given two $d$-intervals~$\maI=\maI _1\wedge\dotsb\wedge \maI _n$ and $\maJ =\maJ _1\wedge\dotsb\wedge \maJ _m$, we define the $d$-interval
\[
\maI\wedge\maJ  = \maI _1\wedge\dotsb\wedge \maI _n\wedge \maJ _1\wedge\dotsb\wedge \maJ _m
\]
and we will regard
$\maI$ and~$\maJ $ as definable subsets of~$\maI \wedge \maJ .$ 

We say that $\maI$ and~$\maJ $ are equal, denoted~$\maI=\maJ $, if $n=m$ and $\maI _i=\maJ _i$ for all~$i=1,\dots , n.$\\

\begin{nrmk}\label{fact interval in d-int}
Let $\maI =\langle I, \langle 0_{\maI}, 1_{\maI}\rangle \rangle $ be a $d$-interval.  If $x,y\in I$ are such that $x\leq _{\maI}y$ then the subset 
$$[x,y]_{\maI}=\{t\in I:x\leq _{\maI}t \leq _{\maI}y\}$$
of the elements of~$\maI$ between $x$ and $y$  is itself a $d$-interval. 

Indeed, let $\maI _i=\langle [a_i, b_i], \langle 0_{\maI _i}, 1_{\maI _i}\rangle \rangle$ be basic $d$-intervals, for $i=1, \ldots , n,$ and suppose that $\maI =\maI _1\wedge \dots \wedge \maI _n.$  Let $i_x, i_y\in \{1,\ldots, n\}$ be such that $x\in [a_{i_x},b_{i_x}]$ and $y\in [a_{i_y},b_{i_y}].$ 
Note that $i_x\leq i_y,$  $x\leq _{\maI}1_{\maI _{i_x}}$ and $0_{\maI _{i_y}}\leq _{\maI}y.$ If  $i_x=i_y$ then
\begin{displaymath}
[x,y]_{\maI} =\left \{ 
\begin{array}{ll}
\langle [x,y], \langle x , y\rangle \rangle  & \textrm{if $a_{i_x}\leq x<y\leq b_{i_x}$}\\
\,\,\, \\
\langle [y, x], \langle x , y\rangle \rangle & \textrm{if $a_{i_x}\leq y<x\leq b_{i_x}$.}
\end{array} \right.
\end{displaymath}

If $i_x<i_y$ then 
\[
[x,y]_{\maI} = [x,1_{\maI _{i_x}}]_{\maI}\wedge {\maI}_{i_x+1} \wedge\dotsb\wedge [0_{\maI _{i_y}},y]_{\maI}.
\]
\end{nrmk}

Below, if $X$ is a locally definable manifold and $Y$ is a definable space, we say that $h:Y\to X$ is a {\it definable continuous map} if for some (equivalently, for every) definable subspace $U$ of $X$ with $h(Y)\subseteq U,$ the map $h:Y\to U$ is a definable continuous map between definable spaces.\\

\begin{defn}
{\em
Let $X$ be a locally definable manifold. A {\it  definable path} $\alpha :\maI \to  X$ is a continuous  (with respect to the topology on $\maI$ given by the order) definable map from some $d$-interval $\maI$ to~$X$. We define $\alpha _0:= \alpha (0_{\maI})$ and $\alpha _1:= \alpha (1_{\maI})$ and call them the endpoints of the definable path $\alpha .$

A definable path $\alpha :\maI \to  X$ is  {\it constant} if  $\alpha _0=\alpha (t)$ for all $t\in \maI$. Below, given a $d$-interval $\maI$ and a point $x\in X,$ we denote by ${\mathbf c}^x_{\maI}$ the constant definable path in $X$ with endpoints $x.$

A definable path $\alpha :\maI \to  X$ is a {\it definable loop} if  $\alpha _0=\alpha _1.$ The {\it inverse} $\alpha ^{-1}$ of a definable path $\alpha :\maI \to  X$ is the definable path  
$$\alpha ^{-1}:=\alpha \circ o_{\maI}^{-1}:\maI \od \to  X.$$ 
A {\it concatenation of two definable paths} $\gamma :\maI \to  X$  and $\delta :\maJ\to  X$ with $\gamma (1_{\maI})=\delta (0_{\maJ})$ is the  definable path $\gamma \cdot \delta :\maI \wedge \maJ \to   X$ with:
\begin{displaymath}
(\gamma \cdot \delta )(t) =\left \{ \begin{array}{ll}
\gamma (t)
& \textrm{if $t\in \maI$}\\
\,\,\, \\
\delta (t)
& \textrm{if $t\in \maJ$.}
\end{array} \right.
\end{displaymath}

We say that $X$ is {\it definably path connected} if for every $u,v$ in $X$ there is a definable path $\alpha :\maI \to   X$ such that $\alpha _0=u$ and $\alpha _1=v.$ \\
}
\end{defn}

In the special case required for our applications we shall prove later, see Corollary \ref{cor pi cover  def jman} (1), that  being definably connected is equivalent to  being definably path connected.

\medskip
Let $X$ be a locally definable manifold and $Y$ a definable space.  Given two definable continuous maps $f,g:Y\to  X$, we say that a definable continuous map $F(t,s):Y\times \maJ\to   X$ is a {\it  definable homotopy between $f$ and $g$} if $f=F_0:=F_{0_{\maJ}}$ and $g=F_1:=F_{1_{\maJ}}$, where $F_s:=F(\cdot ,s)$  for all $s\in \maJ.$ In this case we say that $f$ and $g$ are {\it definably homotopic}, denoted $f\sim g$.\\

\begin{fact}\cite[Remarks 2.8, 2.9 and 2.10]{Edal17}\label{fact hom and conc}
 The definable homotopy $\sim $ is an equivalence relation compatible with concatenation i.e., if $\gamma _i:{\maI } \to X$ and $\delta _i:{\maJ} \to X$ ($i=1,2$) are definable paths with $(\gamma _i)_1=(\delta _i)_0$ for $i=1,2$ and $\gamma _1\sim \gamma _2$ and $\delta _1\sim \delta _2,$ then $\gamma _1\cdot  \delta _1\sim \gamma _2\cdot \delta _2.$
 
 Moreover, if $\gamma : \maI \to X$ is a definable path and $\maJ$ is any $d$-interval, then
\[
{\mathbf c}_{\maI\wedge\maJ }^{\gamma_0}\cdot \gamma \sim \gamma\cdot {\mathbf c}_{\maJ \wedge\maI}^{\gamma_1}.
\]
\end{fact}

Since definable paths need not have the same domain, the notion of homotopic definable paths is not contained in the notion of homotopic definable maps just defined:

\begin{defn}
{\em
Two definable paths $\gamma :\maI \to   X$, $\delta :\maJ\to   X$, with $\gamma _0=\delta _0 $ and $\gamma _1=\delta _1$, are called {\it definably homotopic}, denoted $\gamma \approx \delta ,$ if there  are $d$-intervals $\maI'$ and~$\maJ '$
such that $\maJ '\wedge\maI=\maJ \wedge\maI'$, and there is a definable homotopy 
$${\mathbf c}_{\maJ '}^{\gamma_0}\cdot \gamma
\sim
\delta\cdot  {\mathbf c}_{\maI'}^{\delta_1}$$
fixing the end points (i.e., they are definably homotopic by a definable homotopy $F:\maK \times \maA \to X,$ where $\maK = \maJ '\wedge\maI=\maJ \wedge\maI', $ such that $F(0_{\maK}, s)=\gamma _0=\delta _0$ and $F(1_{\maK}, s)=\gamma _1=\delta _1$ for all $s\in \maA.$)
}
\end{defn}

\medskip

We have:

\begin{nrmk}\cite[Remark 2.11]{Edal17}\label{nrmk sim and approx}
{\em
Let $X$ be a locally definable manifold. If $\delta _i:\maJ \to X$ ($i=1,2$)  are  definable paths such that $\delta _1\sim \delta  _2 ,$ then $\delta _1 \approx \delta _2.$ \\
}
\end{nrmk}

We also have:

\begin{fact}\cite[Proposition 2.13]{Edal17}\label{th-weak-homotopy-eq}
Let $X$ be a locally definable manifold and $x_0, x_1\in X$. Let ${\mathbb P}(X, x_0, x_1)$ denote the set of all  definable paths in $X$ that start at $x_0$ and end at $x_1.$ Then  the restriction of  $\approx$  to ${\mathbb P}(X, x_0, x_1)\times {\mathbb P}(X, x_0, x_1)$ is an equivalence relation on ${\mathbb P}(X, x_0, x_1)$. Moreover, if $\gamma, \gamma', \delta$ and $\delta'$ are definable paths  in  $X$ such that $\gamma _1=\delta _0,$  $\gamma '_1=\delta '_0,$ $\gamma\approx\gamma'$ and $\delta\approx\delta'$, then $\gamma \cdot \delta \approx \gamma' \cdot  \delta' .$
\end{fact}

By \cite[Lemmas 2.14 e 2.15]{Edal17} we have:

\begin{defn}\label{defn fund grp}
{\em 
Let $X$ be a locally definable manifold  and $e_X\in X.$ If ${\mathbb L}(X, e_X)$ denotes the set of all definable loops that start and end at a fixed  element $e_X$ of $X$ (i.e. ${\mathbb L}(X,e_X)={\mathbb P}(X,e_X,e_X)$). 
We define the {\it o-minimal fundamental group} $\pi _1(X, e_X)$ of $X$ by
$$\pi _1(X, e_X):=
\raisebox{1ex}{$\mathsurround=0pt\displaystyle {\mathbb L}(X,e _X)$}
\Big/ \raisebox{-1ex}{$\mathsurround=0pt\displaystyle \approx $}
$$
 with group operation given by $[\gamma ][\delta ]=[\gamma \cdot \delta ],$ the inverse given by $[\gamma ]^{-1}=[\gamma ^{-1}]$ and identity the class of a constant loop at $e_X.$ 

If $f:X\to  Y$ is a locally definable continuous map between two locally definable manifolds with $e_X\in X$ and $e_Y\in Y$ such that $f(e_X)=e_Y$, then we have an induced homomorphism $f_*:\pi _1(X, e_X)\to  \pi _1(Y, e_Y):[\sigma ]\mapsto [f\circ \sigma ]$ with the usual functorial properties. \\
}
\end{defn}

\begin{fact}\cite[Corollary 2.18]{Edal17}\label{fact pi1 and x and connected}
Let $X$ and $Y$ be  locally definable manifolds with $e_X\in X$ and $e_Y\in Y$. Then
\begin{enumerate}
\item
If $X$ is definably path connected then $\pi _1(X, e_X)\simeq \pi _1(X, x)$ for every $x\in X$.
\item
$\pi _1(X, e_X)\times \pi _1(Y, e_Y)\simeq  \pi _1(X\times Y, (e_X,e_Y)).$
\end{enumerate}
\end{fact}



As usual for a definably path connected locally definable manifold $X$ if there is no need to mention a base point $e_X\in X$, then by Fact \ref{fact pi1 and x and connected} (1), we may denote $\pi _1(X,e_X)$ by $\pi _1(X)$.\\

\begin{defn}\label{defn fund grpoid}
{\em 
Let $X$ be a locally definable manifold.  We define the {\it o-minimal fundamental groupoid} $\Pi _1(X)$ of $X$ to be  the small category $\Pi _1(X)$  given by

\begin{displaymath}
\begin{array}{l}
\ob (\Pi _1(X))  = X, \\
\Ho _{\Pi _1(X)}(x_0, x_1)  = \raisebox{1ex}{$\mathsurround=0pt\displaystyle  {\mathbb P}(X, x_0, x_1)$}
\Big/ \raisebox{-1ex}{$\mathsurround=0pt\displaystyle \approx $}\\
\end{array}
\end{displaymath}
We  set $[\gamma ]:=$ the class of $\gamma \in {\mathbb P}(X, x_0,x_1)$. By Fact \ref{th-weak-homotopy-eq}, the small category $\Pi _1(X)$ is indeed a groupoid with operations 
$$\Ho_{\Pi _1(X)}(x_0,x_1)\times \Ho _{\Pi _1(X)}(x_1,x_2)\to \Ho _{\Pi _1(X)}(x_0,x_2)$$ 
given by $[\delta ] \circ [\gamma]=[\gamma \cdot \delta ]$.



If $f:X\to Y$ is a locally definable continuous map between locally definable manifolds, then we have an induced functor $f_*:\Pi _1(X)\to \Pi _1(Y)$ which is a morphism of groupoids sending the object $x\in X$ to the object $f(x)\in Y$ and a morphism $[\gamma ]$ of $\Pi _1(X)$  to the morphism $[f\circ \gamma ]$ of $\Pi _1(Y)$. \\
}
\end{defn}

\end{subsection}

\end{section}

\begin{section}{Main results}\label{sec main}

\begin{subsection}{The main properties}\label{subsec main prop}
The goal in this subsection is to prove the following properties:\\
\begin{propt}\label{propt main new paths}
Let ${\bf P}$ be  the full subcategory of locally definable spaces in ${\mathbb M}$ whose objects are the locally definable  manifolds. Then in the category ${\bf P}$ the following hold:

\begin{itemize}
\item[(P1)]
\begin{itemize}
\item[(a)]
every object of ${\bf P}$ which is definably connected is uniformly definably path connected;
\item[(b)]
given a locally definable covering map $p_X:X\to S$ in ${\bf P}$ then: (i) every definable path $\gamma $ in $S$ has a unique lifting $\tilde{\gamma }$ which is a definable path in $X$ with a given base point; (ii) every definable homotopy $F$ between definable paths $\gamma $ and $\sigma $ in $S$  has a unique lifting $\tilde{F}$ which is a definable homotopy between the definable paths $\tilde{\gamma }$ and $\tilde{\sigma }$ in $X$.
\end{itemize}
\item[(P2)]
Every object of ${\bf P}$ has  admissible covers by definably simply connected, open definable subsets refining any admissible cover by open definable subsets.
\end{itemize} 
\end{propt}

It follows, as observed in the concluding remarks (Section 5) of the paper \cite{eep}, that with (P1) and (P2) above one proves in exactly the same way all  the main results of the paper \cite{eep}, now in the more general context of arbitrary o-minimal structures with definable Skolem functions. These results include all those  mentioned in the Introduction.

In fact,  besides (P1) and (P2) (and their consequences) everything else that is required is, on the one hand, results from \cite{ejp2}, 
which hold in arbitrary o-minimal structures (and for locally definable spaces as well), and  on the other hand,  \cite[Chapter 6, (3.6)]{vdd}, which is used to notice that the domains of the ``good'' definable paths are definably normal. In our case here the domains of the definable paths are Hausdorff, definably compact  definable spaces (Fact \ref{fact d-int space}), definable in the o-minimal structure ${\mathbb M}$ with definable Skolem functions  and so they are definably normal by \cite[Theorem 2.11]{emp}.

The fact that (P1) and (P2) are the only requirements needed to develop the theory presented in \cite{eep} is somewhat natural. Indeed in topology, where we have  good notions of paths and homotopies with the lifting of paths and homotopies property, all one  needs is existence of such nice open covers as in (P2). In the o-minimal context (here and in \cite{eep}),  the role that (P1) (b) and (P2) play is  similar to the role that the analogue properties play in topology.  However, (P2) is often used in combination with the results from \cite{ejp2} mentioned above  to get local definability. Also (P1) (a) is required essentially only once and  to get local definability (see \cite[Proposition 2.18]{eep}), the other places where it is used, it is used to replace definably connected by definably path connected. \\

Below let $\pi :M^{n+1}\to M^n$ be the projection onto the first $n$ coordinates and let $\tau :M^{n+1}\to M: (x,y)\mapsto y$ be the projection onto the last coordinate.\\

We start by proving the existence of continuous definable sections for the projection of an  open  definable subset, over a finite cover by open definable subsets (Theorem \ref{thm continuous-choice} below). But first we recall a few facts.\\ 

The following  is obtained from the definition of cells (\cite[Chapter 3, \S2]{vdd}):\\

\begin{nrmk}\label{nrmk cells perm}
{\em
Let $C\subseteq M^n$ be a $d$-dimensional cell. Then by the definition of cells,  $C$ is a $(i_1, \ldots , i_n)$-cell for some unique sequence $(i_1, \ldots , i_n)$ of $0$'s and $1$'s.  Moreover, if     $\lambda (1)< \dots <\lambda (d)$  are the indices $\lambda \in \{1, \ldots , n\}$ for which $i_{\lambda }=1$ and 
$$p_{\lambda (1), \dots , \lambda (d)}:M^n\to M^d:(x_1, \ldots , x_n)\mapsto (x_{\lambda (1)}, \ldots , x_{\lambda (d)})$$ 
is the projection, then $C':=p_{\lambda (1), \dots , \lambda (d)}(C)$ is an open $d$-dimensional cell in $M^d$ and the restriction $p_C:=p_{\lambda (1), \dots , \lambda (d)|C}:C\to C'$ is a definable homeomorphism  (\cite[Chapter 3, (2.7)]{vdd}). 

Let $\tau (1)<\dots <\tau (n-d)$ be the indices $\tau \in \{1, \ldots , n\}$ for which $i_{\tau }=0.$ For each such $\tau $, by the definition of cells, there is a definable continuous function $h_{\tau }: \pi _{\tau -1}(C)\subseteq M^{\tau -1}\to M$ where, for each $k=1, \ldots , n$,  $\pi _{k}:M^n\to M^{k}$ is the projection onto the first $k$-coordinates. Moreover we have $\pi _{\tau }(C)=\{(x,h_{\tau }(x)): x\in \pi _{\tau -1}(C)\}.$

Let $f=(f_1, \ldots , f_{n-d}):C'\to M^{n-d}$ be the definable continuous map where for each $l=1, \ldots , n-d$ we set $f_l=h_{\tau (l)}\circ \pi _{\tau (l) -1}\circ p_C^{-1}.$ Let  $\sigma :M^n\to M^n:(x_1,\ldots , x_n)\mapsto (x_{\lambda (1)}, \ldots , x_{\lambda (d)}, x_{\tau (1)}, \ldots , x_{\tau (n-d)}).$ Then we clearly have 
$$\sigma (C)=\left\{\left(x,f(x)\right) : x\in C'\right\}.$$
}
\end{nrmk}

\medskip
Recall also the following fact: 

\begin{fact}\cite[Theorem 2.2]{Edal17}\label{thm open cells}
Let $U$ be an open definable subset of $M^n$. 
Then $U$ is a finite union of open definable sets definably homeomorphic, by reordering of coordinates,  to open cells. \\
\end{fact}

\begin{thm}\label{thm continuous-choice}
Let $O$ be an open definable subset of $M^{n+1}.$  Then there is a finite cover $\{U_i:i=1,\ldots, m\}$ of $\pi (O)$ by open definable subsets such that for each $i$ there is a continuous definable section $s_i:U_i\to O$  of $\pi $ (i.e.  $\pi\circ s_i={\rm id}_{U_i}$).
\end{thm}



\pf
By definable Skolem functions, let $s:\pi(O)\to O$ be a definable section, possibly discontinuous, of~$\pi .$ So $s={\rm id}\times t$ for some definable map $t:\pi (O)\to M.$ 

By the cell decomposition theorem (\cite[Chapter 3, (2.11)]{vdd}),  let ${\mathcal C}$ be a cell decomposition of $\pi (O)$ such that $s_{|C}$ is continuous for each cell $C\in {\mathcal C}.$ Thus it is enough to show, by induction on $d,$ that for any cell $C$ of dimension $d$  there are open definable subsets $U_1,\ldots, U_{m}$ of~$\pi(O)$ such that  $C\subseteq \bigcup \{U_i:i=1,\ldots , m\}$ and for each $i$ there  is a continuous definable map $s_i:U_i\to O$   such that $\pi\circ s_i={\rm id}_{U_i}.$  

The result for a zero dimensional cell ($d=0$)   is immediate. For the inductive step, by Remark \ref{nrmk cells perm}, after a reordering of coordinates, we may assume that our cell~$C$ is of the form
\[
C=\left\{\left(x,f(x)\right) : x\in C'\right\}
\]
where $C'$ is a $d$-dimensional open cell in~$M^d$, and $f=(f_1, \ldots , f_{n-d}):C'\to M^{n-d}$ is continuous and definable. 

By definable Skolem functions and the fact that $O$ is open, for each $i=1,\ldots, n-d,$ there  are  definable functions $h_i,g_i:C'\to M $ such that 
 for all $x\in C'$ and for all $y=(y_1,\ldots, y_{n-d})\in M^{n-d}$  the following hold:
\begin{itemize}
\item[-]
for all $i$ we have $h_i(x)<f_i(x)<g_i(x);$
\item[-]
if for all $i$ it holds that $h_i(x)\leq y_i\leq g_i(x)$, then $(x,y, t(x,f(x)))\in O.$
\end{itemize}

Let $\rho:M^n\to M^d$ denote the projection onto the first~$d$ coordinates. By Remark \ref{nrmk cells perm}, the restriction $\rho _{|C}:C\to C'$ to $C$ is an homeomorphism onto $C'$ whose inverse is ${\rm id}\times f:C'\to C.$ 

Let $$K_C=\{(x,y)\in M^{n}: x\in C'\,\,\textrm{and for all}\,\, i,\,\, h_i(x)\leq y_i\leq g_i(x)\}.$$

Then 
\[
s' := {\rm id}\times \left( t\circ \left( {\rm id}\times f \right)\circ \rho \right)
\]
is a continuous definable section on $K_C.$

Let 
$$S=\{x\in C': \,\,\textrm{for all}\,\,i,\,\,\textrm{both}\,\,h_i\,\,\textrm{and}\,\,g_i\,\,\textrm{are continuous at}\,\,x\}.$$ 
Then $S$ is an open definable subset of $C'$ and $C'\setminus S$ has dimension smaller than $d.$ Since $({\rm id}\times f)(S)\subseteq C,$ $C\setminus ({\rm id}\times f)(S)= ({\rm id}\times f)(C'\setminus S)$ and $\dim C=\dim \left(\bar{C}\cap \pi (O)\right),$  it follows that $\left(\bar{C}\cap \pi (O)\right)\setminus ({\rm id}\times f)(S)$  has dimension smaller than~$d.$ (This follows all from properties of o-minimal dimension, see \cite[Chapter 4, \S1]{vdd}).

By the induction hypothesis, $\left(\bar{C}\cap \pi (O)\right)\setminus ({\rm id}\times f)(S)$  can be covered by finitely many open definable subsets of $\pi (O)$  satisfying our requirements. Let $W$ denote the union of these open definable subsets. 

We still have to cover $C\setminus W.$ If $\dim (C\setminus W)<d$ then we apply the induction hypothesis and we are done. So suppose that $\dim (C\setminus W)=d=\dim C.$ Since $\rho _{|C}:C\to C'$ is a definable bijection, we also have $\dim \rho (C\setminus W)=\dim C'.$

Observe that: 

\begin{clm}\label{clm choice1a}
$\rho (C\setminus W)\subseteq S.$
\end{clm}

\pf
Since $({\rm id}\times f)(C'\setminus S)\subseteq W\cap C$, have that $\rho(W\cap C)$  contains $C'\setminus S.$  Since $\rho _{|C}$ is a bijection, $\rho (C\setminus W)=C'\setminus \rho (W\cap C)\subseteq S.$  \qed \\

Let $$K_{C\setminus W}=\{(x,y)\in M^{n}: x\in \rho (C\setminus W)\,\,\textrm{and for all}\,\, i,\,\, h_i(x)\leq y_i\leq g_i(x)\}.$$
Let $V$ be the interior of $\rho (C\setminus W)$ and let 
$$U=\{(x,y)\in M^{n}: x\in V\,\,\textrm{and for all}\,\, i,\,\, h_i(x)< y_i< g_i(x)\}.$$


\begin{clm}\label{clm choice2a}
Then $V\neq \emptyset ,$ $U$ is open and $ U\subseteq K_{C\setminus W}$.
\end{clm}

\pf
The fact that $V\neq \emptyset $ follows from the fact that $C'$ is open and $\dim \rho (C\setminus W)=\dim C'$ (see \cite[Chapter 4, (1.9)]{vdd}). Since $V\subseteq \rho (C\setminus W),$  we have $U\subseteq K_{C\setminus W}.$ Since $\rho (C\setminus W)\subseteq S$ (Claim \ref{clm choice1a}), for all $i$, $h_i, g_i:V\to M$ are continuous definable maps. Therefore,  $U$ is open in $\rho ^{-1}(V)$ and so it is open since $V$ is open in $C'$ and $C'$ is open in $M^d.$ \qed \\







Since  $U \subseteq K_{C\setminus W}\subseteq K_C,$ we have that $s'$ is continuous on $U.$ We still have to cover $(C\setminus W)\setminus U.$ However since $C'$ is open and $\dim \rho (C\setminus W)=\dim C'$ it follows that $\dim (\rho (C\setminus W)\setminus V)<d$ (see \cite[Chapter 4, (1.9)]{vdd}). So $\dim ((C\setminus W)\setminus U)<d$  and we conclude by the induction hypothesis.
\qed \\

Let $Z\subseteq M^{n+1}$ be a definable set, $X$ a  definable space and let $f:X\to  \pi (Z)$ be a definable  map.  A {\it lifting of $f$} is a continuous map $g:X\to  Z$  such that $\pi \circ g=f$. \\

\begin{lem}\label{lem unique hom lifting}
Let $Z\subseteq M^{n+1}$ be a definable set.  Let $\maI =\langle I, \langle 0_{\maI},1_{\maI}\rangle \rangle $ be a basic $d$-interval, $\gamma :\maI  \to \pi (Z)$ be a definable path in $\pi (Z)$ and let $\alpha , \beta :{\maI}\to Z$ be definable paths lifting $\gamma .$ 
Suppose that for all $t\in I$ we have $\{(\gamma (t), t):\min\, (\tau \circ \alpha (t),\tau\circ\beta (t))\leq t\leq \max\, (\tau \circ \alpha (t),\tau\circ\beta (t))\}\subseteq Z.$
Then $\alpha \sim \beta .$ Moreover, if $\alpha _0=\beta _0$ and $\alpha _1=\beta _1,$  then the definable homotopy $\alpha \sim \beta $ fixes the endpoints.
\end{lem}

\pf
Let $\mu :\maI \to C$ be the definable path given by 
$$\mu (t)= (\gamma (t), \min\,(\tau \circ \alpha (t), \tau\circ \beta (t))).$$
Let $a=\max \{\max \,(\tau\circ \alpha (t), \tau \circ \beta (t)):t\in I\}$ and let $b=\min \{\min \,(\tau\circ \alpha (t), \tau \circ \beta (t)):t\in I\}.$  If $a=b$ then $\alpha =\beta ,$ so  we may assume that $b<a.$ Consider the basic $d$-interval $\maK=\langle [b,a], \langle a,b\rangle \rangle .$ Let $F:\maI\times \maK \to C$ be the continuous definable map given by 
$$F(t,r)=(\gamma (t), \max \,(\tau \circ \mu (t), \min \, (\tau \circ \alpha (t),r))).$$
Then $F_0=\alpha $ and $F_1=\mu .$ Therefore, $\alpha \sim \mu .$ Similarly, $\beta \sim \mu .$ Hence $\alpha \sim \beta$.
\qed \\

We will want to show that a given definable path is definably homotopic to a second definable path, but we cannot do it because their domains are not the same. For this we will modify the first definable path by  patching to it  appropriate constants that allow it to ``wait'' for the second definable path. That ``waiting period'' can occur either before or after (or both) appropriate basic $d$-intervals that composes the $d$-interval of the domain of the first path.

We say that a definable path $\gamma ':\maI '\to X$ is obtained from a definable path $\gamma :\maI\to X$ by {\it modifying with constants} if there are $0_{\maI}=t_0<_{\maI}t_1<_{\maI}\cdots <_{\maI}t_r=1_{\maI}$ and there are $d$-intervals ${\maJ}_0,\ldots , {\maJ}_r$ such that,  if  $\gamma^i:=\gamma _{|[t_i,t_{i+1}]_{\maI}}$ for $i=0,1,\ldots r-1,$ then  one of the following three cases holds:
\begin{enumerate}
\item
 $\maI ' ={\maJ}_0\wedge [t_0,t_1]_{\maI}\wedge \ldots \wedge {\maJ}_{r-1}\wedge [t_{r-1}, t_r]_{\maI}$ and $\gamma '= (\mathbf{c}^{\gamma ^0 _0}_{{\maJ}_0}\cdot \gamma ^0)\cdot \cdots   \cdot (\mathbf{c}^{\gamma ^{r-1}_0}_{{\maJ}_{r-1}}\cdot \gamma ^{r-1}).$
 \item
  $\maI ' ={\maJ}_0\wedge [t_0,t_1]_{\maI}\wedge \ldots \wedge {\maJ}_{r-1}\wedge [t_{r-1}, t_r]_{\maI}\wedge {\maJ}_r$ and $\gamma '= (\mathbf{c}^{\gamma ^0 _0}_{{\maJ}_0}\cdot \gamma ^0)\cdot \cdots   \cdot (\mathbf{c}^{\gamma ^{r-1}_0}_{{\maJ}_{r-1}}\cdot \gamma ^{r-1})\cdot \mathbf{c}^{\gamma ^{r-1}_1}_{{\maJ}_{r}}.$
 \item
  $\maI ' =[t_0,t_1]_{\maI}\wedge {\maJ}_1\wedge [t_1,t_2]_{\maI}\wedge\ldots \wedge {\maJ}_{r-1}\wedge [t_{r-1}, t_r]_{\maI}\wedge {\maJ}_r$ and $\gamma '= \gamma ^0\cdot (\mathbf{c}^{\gamma ^1 _0}_{{\maJ}_1}\cdot \gamma _1)\cdot \cdots   \cdot (\mathbf{c}^{\gamma ^{r-1}_0}_{{\maJ}_{r-1}}\cdot \gamma ^{r-1})\cdot \mathbf{c}^{\gamma ^{r-1}_1}_{{\maJ}_{r}}.$ \\
 \end{enumerate}

\begin{lem}\label{lem mod by const}
If a definable path $\gamma ':\maI '\to X$ is obtained from a definable path $\gamma :\maI\to X$ by  modifying with constants, then $\gamma \approx \gamma '.$ 
\end{lem}

\pf
Consider the first case. By Fact \ref{fact hom and conc}, 
$$\mathbf{c}^{(\mathbf{c}^{\gamma ^{i}_0}_{{\maJ}_{i}}\cdot \gamma ^{i})_0}_{[t_{i}, t_{i+1}]_{\maI}}\cdot (\mathbf{c}^{\gamma ^{i}_0}_{{\maJ}_{i}}\cdot \gamma ^{i})=\mathbf{c}_{[t_{i}, t_{i+1}]_{\maI}\wedge {\maJ}_{i}}\cdot \gamma ^{i}
\sim \gamma ^i\cdot \mathbf{c}^{\gamma ^i_1}_{{\maJ}_{i}\wedge [t_i,t_{i+1}]_{\maI}}$$
and so $(\mathbf{c}^{\gamma ^{i}_0}_{{\maJ}_{i}}\cdot \gamma ^{i})\approx \gamma ^i.$ Since $\gamma =\gamma ^0\cdot \gamma ^1\cdot \cdots \cdot \gamma ^{r-1},$   we conclude by Fact \ref{th-weak-homotopy-eq} that $\gamma \approx \gamma '.$

The third case is similar and the second case follows from the first two cases by transitivity of $\approx.$
\qed \\

\begin{nrmk}\label{remark interval def normal}
{\em 
Let $I=\Pi _{i=1}^n[t_{i-1}, t_i] \subseteq M^n$ be a product of intervals. By \cite[Theorem 2.1]{ps}, $I$ is a Hausdorff, definably compact  definable space. Since ${\mathbb M}$ has definable Skolem functions, it follows that $I$ is  definably normal (\cite[Theorem 2.11]{emp}).
}
\end{nrmk}

\begin{lem}\label{lem lifting paths}
Let $C\subseteq M^{n+1}$ be an open cell.  Then the following hold.
\begin{enumerate}
\item
Let $\gamma :\maI  \to \pi (C)$ be a definable path in $\pi (C).$ Let $(x,y)\in C$ be such that $x=\gamma _{0}.$ Then there is a  definable path $\gamma ':{\maI}'\to \pi (C)$ obtained from $\gamma $ by modifying with constants and there is a lifting $\beta:\maI '\to C$  of $\gamma '$ with $\beta _{0}=(x,y)$.

\item
Suppose that $F: \maI \times \maJ \to \pi (C)$ is a definable homotopy  between the definable paths $\gamma ,\sigma :{\maI}\to \pi (C)$ in $\pi (C)$. Then there are definable paths $\gamma ',\sigma ':{\maI}'\to \pi (C)$ obtained by modifying with constants $\gamma $ and $\sigma $ respectively, and there are liftings $\beta , \tau :{\maI}'\to C$ of $\gamma '$ and $\sigma '$ respectively such that $\beta \sim \tau .$ 
\end{enumerate}
\end{lem}

\pf
By Theorem \ref{thm continuous-choice}, there is a finite cover $\{U_i:i=1,\ldots, m\}$ of $\pi (C)$ by open definable subsets such that for each $i$ there is a continuous definable section $s_i:U_i\to C$  of the projection $\pi $ (i.e.  $\pi\circ s_i={\rm id}_{U_i}$).\\

(1) First we assume that $\maI$ is a basic $d$-interval $\langle [a,b], \langle 0_{\maI}, 1_{\maI}\rangle \rangle .$ We may also assume that the definable total order  $<_{\maI}$ on the domain $[a, b]$ of $\maI$ is $<.$ If not, the argument is similar, one just has to construct the lifting from right to left instead of from left to right.

Let $L=\{l: \gamma ([a,b])\cap U_l\neq \emptyset \}.$ Then $[a,b]\subseteq  \bigcup_{l\in L}\gamma^{-1}(U_l)$, with  the $\gamma^{-1}(U_l)$'s open in $[a,b]$. So by Remark \ref{remark interval def normal} and the shrinking lemma (\cite[Corollary 2.12]{emp},  \cite[Chapter 6, (3.6)]{vdd}),  for each $l\in L$ there is $W_l\subset [a,b]$, open in $[a,b]$ such that $W_l\subset\overline{W_l}\subset\gamma^{-1}(U_l)$ and $[a,b]\subseteq  \bigcup_{l\in L} W_l$. Therefore,  there are $a=t_0<t_1<\cdots <t_r=b$ such that for each $i=0,\dots,r-1$ we have $\gamma([t_i,t_{i+1}]) \subset U_{l(i)}$ (and $\gamma(t_{i+1})\in U_{l(i)}\cap U_{l(i+1)}$).  For each $i=0,\dotsc, r-1$ let $\gamma^i:=\gamma _{|[t_i,t_{i+1}]}.$ Then $s_{l(i)}\circ\gamma^i :[t_i, t_{i+1}]\to C$ is a definable path joining $s_{l(i)}(\gamma^i_0)$ to~$s_{l(i)}(\gamma^i_1)$. Now, observe that the points $s_{l(i)}(\gamma^i_1)$ and~$s_{l(i+1)}(\gamma^{i+1}_0)$ differ just by the $y$~coordinate, and the same happens with the points $(x,y)$ and $s_{l(0)}(\gamma^0_0).$ Such pairs of points are clearly connected by (definable) vertical paths. Let $\nu ^0:{\maJ}_0\to C$ be the vertical path with $\nu ^0_0=(x,y)$ and $\nu ^0_1=s_{l(0)}(\gamma ^0_0)$ and for $i=1,\ldots r-1$ let $\nu ^i:{\maJ}_i\to C$ the vertical path with $\nu ^i_0=s_{l(i-1)}(\gamma ^{i-1}_1)$ and $\nu ^i_1=s_{l(i)}(\gamma ^{i}_0).$ Let ${\maI}'={\maJ}_0\wedge [t_0,t_1]\wedge \ldots \wedge {\maJ}_{r-1}\wedge [t_{r-1},t_r]$ and let $\beta =\nu ^0\cdot (s_{l(0)}\circ \gamma ^0)\cdot \cdots   \cdot \nu ^{r-1}\cdot (s_{l(r-1)}\circ \gamma ^{r-1}).$ Let $\gamma ':{\maI}'\to \pi (C)$ be given by $\gamma '= (\mathbf{c}^{\gamma ^0 _0}_{{\maJ}_0}\cdot \gamma ^0)\cdot \cdots   \cdot (\mathbf{c}^{\gamma ^{r-1}_0}_{{\maJ}_{r-1}}\cdot \gamma ^{r-1}).$ Then the definable path $\beta :\maI ' \to C$ in $C$ is a lifting of $\gamma '$ such that $\beta  _{0}=(x,y)$.

Now if $\maI =\maI _1\wedge \ldots \wedge \maI _k$ with each $\maI_i$ a basic $d$-interval apply the previous process to  $\gamma_{|\maI _1}$ to get $\beta ^1,$ $\gamma '^1$ with  $\beta ^1(0_{{\maI}_1})=(x,y)$ and repeat the process for each $\gamma _{|\maI _{i+1}}$ with $\beta ^i(1_{\maI_i})$ instead of $(x,y)$. Patch these together to obtain $\beta $ and $\gamma '.$ \\

(2) Here we cannot just apply (1) to $\gamma $ and $\sigma $ since the liftings $\beta $ and $\tau $ of the corresponding modifications by constants need not be definably homotopic. Instead we need to use the definable homotopy $F$ to build a sequence of definable paths which are ``close'' enough to guarantee that the liftings of the corresponding modifications by constants are definably homotopic. The sequence of definable paths will be obtained from a cell decomposition of the domain of the definable homotopy $F$ compatible with the pull backs of the open definable subsets of $\pi (C)$ on which  we have continuous definable sections (Theorem \ref{thm continuous-choice}) and the liftings of modifications by constants of these definable  paths will be obtained using the continuous definable sections.

First assume  that $\maJ$ is a basic $d$-interval $\langle [c,d], \langle 0_{\maJ}, 1_{\maJ}\rangle \rangle .$ We may also assume that the definable total order  $<_{\maJ}$ on the domain $[c, d]$ of $\maJ$ is $<.$ If not, the argument is similar, one just has to construct the lifting from top to bottom  instead of from bottom  to top.

To proceed we also assume that $\maI$ is a basic $d$-interval $\langle [a,b], \langle 0_{\maI}, 1_{\maI}\rangle \rangle .$ We may furthermore  assume that the definable total order  $<_{\maI}$ on the domain $[a, b]$ of $\maI$ is $<.$ If not the argument is similar, one just has to construct the lifting from right to left instead of from left to right. 

Let $L=\{l:F([a,b]\times [c,d])\cap U_l\neq \emptyset \}.$ Then $[a,b]\times [c,d] \subseteq \bigcup_{l\in L}F^{-1}(U_l)$, with  the $F^{-1}(U_l)$'s open in $[a, b]\times [c,d]$. So by Remark \ref{remark interval def normal} and the shrinking lemma (\cite[Corollary 2.12]{emp}, \cite[Chapter 6, (3.6)]{vdd}),   we have that  for each $l\in L$ there is $W_l\subset [a, b]\times [c,d]$, open in $[a, b]\times [a,d]$ such that $W_l\subset\overline{W_l}\subset F^{-1}(U_l)$ and $[a, b]\times[c, d] \subseteq \bigcup_{l\in L}W_l$. Now take a cell decomposition of  $[a,b]\times [c,d]$ compatible with the $W_l$'s. This cell decomposition is given by  a decomposition  $a=t_0<t_1<\cdots <t_r=b$ of $[a, b]$ together with definable continuous functions $f_{i,j}:[t_i, t_{i+1}] \to [c,d]$ for $i=0,\ldots , r-1$ and $j=0,\ldots , k_i$ such that: (i) $f_{i,0}<f_{i,1}< \ldots <f_{i, k_i}$ for $i=0,\ldots , r-1;$ (ii) $f_{i,0}=c$ and $f_{i, k_i}=d$ for $i=0,\ldots , r-1;$ (iii) the two-dimensional cells are of the form $C_{i,j=}(f_{i,j}, f_{i, j+1})_{(t_i, t_{i+1})}.$ For each two-dimensional cell $C_{i,j}$ and each $l(i,j)$ such that $C_{i,j}\subset W_{l(i,j)}$, we have $F(\overline{C_{i,j}})\subset U_{l(i,j)}$ and  for any two-dimensional cells $C_{i,j}$ and $C_{i',j'}$ in $[a, b]\times [c,d]$,   and for each $l(i,j), l(i',j'),$ such that $C_{i,j}\subset W_{l(i,j)}$ and $C_{i',j'}\subset W_{l(i',j')}$ we also have  $F(\overline{C_{i,j}}\cap\overline{C_{i',j'}})\subset U_{l(i,j)}\cap U_{l(i',j')}$. 

Fix a sequence $\bar{j}=(j_0, j_1,\ldots, ,j_{r-1})$ with $j_i\in \{0,1,\ldots, n_i\}$ for each $i=0,\ldots, r-1.$  For each $i=1, \ldots , r$ let  $r^{\bar{j},i}:\maR _{\bar{j},i}\to [a,b]\times [c,d]$ be the vertical path in $[a,b]\times [c,d]$ with $r^{\bar{j},i}_0=(t_i,f_{i-1,j_{i-1}}(t_i))$ and $r^{\bar{j},i}_1=(t_i, f_{i-1,\min\,(j_{i-1}+1,n_{i-1})}(t_i)).$ This is  the vertical path, going up on the right hand side of the ``cell'' 
$$[f_{i-1,j_{i-1}}, f_{i-1,\min\,(j_{i-1}+1,n_{i-1})}]_{(t_{i-1},t_i)}.$$

For each $i=0, \ldots , r-1$ let  $l^{\bar{j},i}:\maL _{\bar{j},i}\to [a,b]\times [c,d]$ be the vertical path in $[a,b]\times [c,d]$ with $l^{\bar{j},i}_0=(t_i,f_{i,j_i}(t_i))$ and $l^{\bar{j},i}_1=(t_i, f_{i,\min\,(j_i+1,n_i)}(t_i)).$ This is  the vertical path, going up on the left hand side of the ``cell'' 
$$[f_{i,j_{i}}, f_{i,\min\,(j_{i}+1,n_{i})}]_{(t_{i},t_{i+1})}.$$

For $i=0,1,\ldots r-1$ let 
$$a^{\bar{j},i}:[t_i,t_{i+1}]\to [a,b]\times [c,d]$$
be the definable path given by
$$a^{\bar{j},i}(t)=(t, f_{i,j_i}(t)).$$

For each $i=1, \ldots , r-1$ let $u^{\bar{j},i}:\maU _{\bar{j},i}\to [a,b]\times [c,d]$ be the vertical path in $[a,b]\times [c,d]$ with $u^{\bar{j},i}_0=(t_i,f_{i-1,j_{i-1}}(t_i))=a^{\bar{j},i-1}_1$ and $u^{\bar{j},i}_1=(t_i, f_{i,j_{i}}(t_i))=a^{\bar{j},i}_0.$  

Let
$$a^{\bar{j}}:[t_0,t_1]\wedge {\maU}_{\bar{j},1}\ldots \wedge {\maU}_{\bar{j},r-1}\wedge [t_{r-1},t_r]\to [a,b]\times [c,d]$$
be the definable path given by
$$a^{\bar{j}}=a^{\bar{j},0}\cdot u^{\bar{j},1}\cdot a^{\bar{j},1}\cdot \cdots \cdot u^{\bar{j},r-1}\cdot a^{\bar{j},r-1}.$$

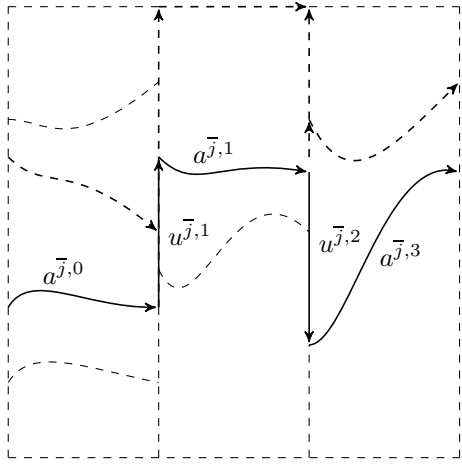
\vspace{1cm}
\begin{center}
\label{pic1}
\begin{tikzpicture}
  \draw[->,>=stealth',shorten >=1pt,auto,node distance=2.8cm,
                    semithick] (0,1) .. controls (0.3,1.5) and (1,1) ..  node[above] {$a^{\bar{j},0}$}  (2,1);
    \draw[->,>=stealth',shorten >=1pt,auto,node distance=2.8cm,
                    semithick] (2,1) .. controls (2,1) and (2,3) .. node[right] {$u^{\bar{j},1}$} (2,3);
    \draw[->,>=stealth',shorten >=1pt,auto,node distance=2.8cm,
                    semithick] (2,3) .. controls (2.5,2.5) and (2.8,3) ..  node[above] {$a^{\bar{j},1}$}  (4,2.8);
    \draw[->,>=stealth',shorten >=1pt,auto,node distance=2.8cm,
                    semithick] (4,2.8) .. controls (4,2) and (4,2) .. node[right] {$u^{\bar{j},2}$} (4,0.5);
     \draw[->,>=stealth',shorten >=1pt,auto,node distance=2.8cm,
                    semithick] (4,0.5) .. controls (4.5,0.5) and (5,3) ..  node[below, right] {$a^{\bar{j},3}$}  (6,2.8);
     \draw[dashed] (0,-1) -- (0,5);
     \draw[dashed] (2,-1) -- (2,1);
     \draw[->,>=stealth',shorten >=1pt,auto,node distance=2.8cm,
                    semithick][dashed] (2,1) -- (2,5);
     \draw[dashed] (4,-1) -- (4,0.5);
     \draw[->,>=stealth',shorten >=1pt,auto,node distance=2.8cm,
                    semithick][dashed] (4,3) -- (4,3.5);
     \draw[->,>=stealth',shorten >=1pt,auto,node distance=2.8cm,
                    semithick][dashed] (4,3.5) -- (4,5);
     \draw[dashed] (6,-1) -- (6,5);
     \draw[dashed] (0,-1) -- (6,-1);
     \draw[dashed] (0,5) -- (2,5);
     \draw[dashed] (4,5) -- (6,5);
     \draw[->,>=stealth',shorten >=1pt,auto,node distance=2.8cm,
                    semithick][dashed] (2,5) -- (4,5);
     \draw[dashed] (0,0) .. controls (0.3,0.5) and (1,0.2) ..  node[above] {}  (2,0);
      \draw[->,>=stealth',shorten >=1pt,auto,node distance=2.8cm,
                    semithick][dashed] (0,3) .. controls (0.5,2.5) and (0.8,3) ..  node[above] {}  (2,2);
      \draw[dashed] (0,3.5) .. controls (0.5,3.5) and (0.8,3) ..  node[above] {}  (2,4);
    \draw[dashed] (2,1.5) .. controls (2.5,0.5) and (3,3) ..  node[above] {}  (4,2);
     \draw[->,>=stealth',shorten >=1pt,auto,node distance=2.8cm,
                    semithick][dashed] (4,3.5) .. controls (4.5,2.5) and (5,3) ..  node[below, right] {}  (6,4);

\end{tikzpicture}
\end{center}
\vspace{1cm}

For $k=0,1,\ldots r-1$  let $\bar{j}[k]$ the sequence which is equal to $\bar{j}$ except in position $k$ where it is $\min\,(j_k+1,n_k).$ \\

By  Lemma  \ref{lem unique hom lifting}, we have:\\

\begin{nrmk}\label{nrmk as}
If  $y$ is a point on the right hand side of $C_{k,j_k}$ let $y^-:{\maY}^-\to [a,b]\times [c,d]$ be the vertical path such that $y^-_0=a^{\bar{j}[k],k}_1$ and $y^-_1=y$ and let $y^+:{\maY}^+\to [a,b]\times [c,d]$ be the vertical path such that $y^+_0=a^{\bar{j},k}_1$ and $y^+_1=y.$ 

If  $x$ is a point on the left hand side of $C_{k,j_k}$ let $x^+:{\maX}^+\to [a,b]\times [c,d]$ be the vertical path such that $x^+_0=x$ and $x^+_1=a^{\bar{j}[k],k}_0$  and let $x^-:{\maX}^-\to [a,b]\times [c,d]$ be the vertical path such that $x^-_0=x$ and $x^-_1=a^{\bar{j},k}_0.$ 
 
Then 

\begin{itemize}
\item[(i)]
$$a^{\bar{j}[0],0}\cdot y^-\cdot \mathbf{c}^{y}_{{\maY}^+}\sim a^{\bar{j},0}\cdot  \mathbf{c}^{a^{\bar{j},0}_1}_{{\maY}^-}\cdot y^+$$
 and  $a^{\bar{j}[0],i}=a^{\bar{j},i}$ for $i\neq 0.$
\item[(ii)]
$$x^+\cdot \mathbf{c}^{x^+_1}_{\maX ^-}\cdot a^{\bar{j}[k],k} \cdot y^- \cdot \mathbf{c}^y_{\maY ^+}\sim \mathbf{c}^{x}_{\maX ^+}  \cdot x^- \cdot a^{\bar{j},k} \cdot \mathbf{c}^{a^{\bar{j},k}_1}_{\maY ^-}\cdot y^+$$
and $a^{\bar{j}[k],i}=a^{\bar{j},i}$ for $i\neq k.$
\item[(iii)]
$$x^+\cdot \mathbf{c}^{x^+_1}_{\maX ^-}\cdot a^{\bar{j}[r-1],r-1} \sim \mathbf{c}^{x}_{\maX ^+}  \cdot x^- \cdot a^{\bar{j},r-1}$$
and $a^{\bar{j}[r-1],i}=a^{\bar{j},i}$ for $i\neq r-1.$
\end{itemize}
\end{nrmk}

On the other hand we also have, see the picture on page \pageref{pic1}:\\

\begin{nrmk}\label{nrmk us}
We have:
\begin{itemize}
\item[(i)]
$u^{\bar{j}[0],1}=(r^{\bar{j},1})^{-1}\cdot u^{\bar{j},1}$ or $u^{\bar{j},1}=r^{\bar{j},1}\cdot u^{\bar{j}[0],1}$ or $u^{\bar{j}[0],1}_1=u^{\bar{j},1}_1$ is on the   right hand side of $C_{0,j_0}.$ \\
Furthermore,  $u^{\bar{j}[0],i}=u^{\bar{j},i}$ for $i\neq 1.$
\item[(ii)]
$u^{\bar{j}[k],k}=u^{\bar{j},k}\cdot l^{\bar{j},k}$ or $u^{\bar{j},k}= u^{\bar{j}[k],k}\cdot (l^{\bar{j},k})^{-1}$ or $u^{\bar{j}[k],k}_0=u^{\bar{j},k}_0$ is on the  left hand side of $C_{k,j_k}.$ \\
$u^{\bar{j}[k],k+1}=(r^{\bar{j},k+1})^{-1}\cdot u^{\bar{j},k+1}$ or $u^{\bar{j},k+1}=r^{\bar{j},k+1}\cdot u^{\bar{j}[k],k+1}$ or $u^{\bar{j}[k],k+1}_1=u^{\bar{j},k+1}_1$ is on the   right hand side of $C_{k,j_k}.$\\
Furthermore, $u^{\bar{j}[k],i}=u^{\bar{j},i}$ for $ i\neq k, k+1.$
\item[(ii)]
$u^{\bar{j}[r-1],r-1}= u^{\bar{j},r-1}\cdot l^{\bar{j},r-1}$ or $u^{\bar{j},r-1}= u^{\bar{j}[r-1],r-1}\cdot (l^{\bar{j},r-1})^{-1}$ or $u^{\bar{j}[r-1],r-1}_0$ $=u^{\bar{j},r-1}_0$ is on the  left hand side of $C_{r-1,j_{r-1}}.$ \\
 Furthermore, $u^{\bar{j}[r-1],i}=u^{\bar{j},i}$ for $i\neq r-1.$\\
\end{itemize}
\end{nrmk}

Let $\mu ^{\bar{j},i}=F\circ u^{\bar{j},i}:\maU _{\bar{j},i}\to \pi (C),$ $\rho ^{\bar{j},i}=F\circ r^{\bar{j},i}:\maR _{\bar{j},i}\to \pi(C),$ $\lambda ^{\bar{j},i}=F\circ l^{\bar{j},i}:\maL _{\bar{j},i}\to \pi(C),$ $\alpha ^{\bar{j},i}=F\circ a^{\bar{j},i}:[t_i,t_{i+1}]\to \pi (C)$ and 
$$\alpha ^{\bar{j}}=F\circ a^{\bar{j}}:[t_0,t_1]\wedge {\maU}_{\bar{j},1}\ldots \wedge {\maU}_{\bar{j},r-1}\wedge [t_{r-1},t_r]\to \pi (C).$$
Then
$$\alpha ^{\bar{j}}=\alpha ^{\bar{j},0}\cdot \mu^{\bar{j},1}\cdot a^{\bar{j},1}\cdot \cdots \cdot \mu^{\bar{j},r-1}\cdot \alpha ^{\bar{j},r-1},$$
$\alpha ^{\bar{0}}=\gamma $ and $\alpha ^{\bar{n}}=\sigma $ where $\bar{0}=(0,0, \ldots, 0)$ and $\bar{n}=(n_0,n_1,\ldots ,n_{r-1}).$
$\,$\\

Applying $F$ to  Remarks \ref{nrmk as} and \ref{nrmk us} we obtain:

\begin{clm}
For every $\bar{j}$ and every $k,$ after modifying with constants, that we ignore for simplicity, 
$$\alpha ^{\bar{j}}\sim \alpha ^{\bar{j}[k]}.$$
Moreover, since for every $\bar{j'}$ there are $k_1, \ldots, k_m$ such that $\bar{j'}=\bar{j}[k_1]\ldots [k_m]$ we also have, by transitivity of $\sim $ (Fact \ref{fact hom and conc}), ignoring modifications by constants,  
$$\alpha ^{\bar{j}}\sim \alpha ^{\bar{j'}}.$$
\end{clm}

\pf
Suppose that $k=0.$ Then $\alpha ^{\bar{j}}=\alpha ^{\bar{j},0}\cdot \mu^{\bar{j},1}\cdot \alpha '$ and $\alpha ^{\bar{j}[0]}=\alpha ^{\bar{j}[0],0}\cdot \mu^{\bar{j}[0],1}\cdot \alpha '.$  By Remark \ref{nrmk us} (i) we have three cases. Suppose that  $u^{\bar{j}[0],1}=(r^{\bar{j},1})^{-1}\cdot u^{\bar{j},1}.$ Then $\mu ^{\bar{j}[0],1}=(\rho ^{\bar{j},1})^{-1}\cdot \mu ^{\bar{j},1}.$ Therefore, $\alpha ^{\bar{j}[0]}=\alpha ^{\bar{j}[0],0}\cdot (\rho ^{\bar{j},1})^{-1}\cdot \mu ^{\bar{j},1}\cdot \alpha '.$  Since $a ^{\bar{j}[0],0}\cdot (r ^{\bar{j},1})^{-1}\sim a^{\bar{j},0}$ by Remark \ref{nrmk as}, we get $\alpha ^{\bar{j}[0],0}\cdot (\rho ^{\bar{j},1})^{-1}\sim \alpha ^{\bar{j},0}$ and the result follows. The other cases are similar.

Suppose that $k=1,\ldots, r-2.$ Then $\alpha ^{\bar{j}}=\alpha ''\cdot \mu ^{\bar{j},k}\cdot \alpha ^{\bar{j},k}\cdot \mu^{\bar{j},k+1}\cdot \alpha '$ and $\alpha ^{\bar{j}[k]}=\alpha ''\cdot \mu ^{\bar{j}[k],k}\cdot \alpha ^{\bar{j}[k],k}\cdot \mu^{\bar{j}[k],k+1}\cdot \alpha '.$ By Remark \ref{nrmk us} (i) we have several cases. Suppose that $u^{\bar{j},k}= u^{\bar{j}[k],k}\cdot (l^{\bar{j},k})^{-1}$ and $u^{\bar{j}[k],k+1}_1=u^{\bar{j},k+1}_1$ is on the   right hand side of $C_{k,j_k}.$ Then $\mu^{\bar{j},k}= \mu^{\bar{j}[k],k}\cdot (\lambda^{\bar{j},k})^{-1}$ and so $\alpha ^{\bar{j}}=\alpha ''\cdot \mu^{\bar{j}[k],k}\cdot (\lambda^{\bar{j},k})^{-1}\cdot \alpha ^{\bar{j},k}\cdot \mu^{\bar{j},k+1}\cdot \alpha '.$ Since by Remark \ref{nrmk as},  
$(l^{\bar{j},k})^{-1}\cdot a ^{\bar{j},k}\cdot u^{\bar{j},k+1} \sim a ^{\bar{j}[k],k}\cdot u^{\bar{j}[k],k+1}$ we obtain $(\lambda^{\bar{j},k})^{-1}\cdot \alpha ^{\bar{j},k}\cdot \mu^{\bar{j},k+1} \sim \alpha ^{\bar{j}[k],k}\cdot \mu^{\bar{j}[k],k+1}$ and the result follows. The other cases are similar.

For $k=r-1$ the argument is the same.
\qed \\

By Lemmas \ref{lem unique hom lifting} and \ref{lem mod by const} and the transitivity of $\sim $ (Fact \ref{fact hom and conc}), to finish the proof it is enough to show that  after modifying with constants, $\alpha ^{\bar{j}}$ and $\alpha ^{\bar{j}[k]}$ have liftings $\beta ^{\bar{j}}$ and $\beta ^{\bar{j}[k]}$ respectively such that
$$\beta ^{\bar{j}}\sim \beta ^{\bar{j}[k]},$$
ignoring modifications by constants.\\

\medskip
\noindent
Suppose that $k=0.$ Then $\alpha ^{\bar{j}}=\alpha ^{\bar{j},0}\cdot \mu^{\bar{j},1}\cdot \alpha '$ and $\alpha ^{\bar{j}[0]}=\alpha ^{\bar{j}[0],0}\cdot \mu^{\bar{j}[0],1}\cdot \alpha '.$ By Remark \ref{nrmk us} (i) we have three cases. Suppose that  $u^{\bar{j}[0],1}=(r^{\bar{j},1})^{-1}\cdot u^{\bar{j},1}.$ Then $\mu ^{\bar{j}[0],1}=(\rho ^{\bar{j},1})^{-1}\cdot \mu ^{\bar{j},1}.$ Therefore, $\alpha ^{\bar{j}[0]}=\alpha ^{\bar{j}[0],0}\cdot (\rho ^{\bar{j},1})^{-1}\cdot \mu ^{\bar{j},1}\cdot \alpha '.$ 

By (1) let $\mu $ be a lifting of (a modification by constants of) $\mu ^{\bar{j},1}$ and let $\beta '$ be a lifting of (a modification by constants of) $\alpha '$ such that $\beta '_0=\mu _1.$ Recall that $F(\bar{C_{0,j_0}})\subseteq U_{l(0,j_0)}.$ Thus $\rho ^{\bar{j},1}=F\circ r^{\bar{j},1}:\maR _{\bar{j},1}\to U_{l(0,j_0)} \subseteq\pi(C),$   $\alpha ^{\bar{j},0}=F\circ a^{\bar{j},0}:[t_0,t_{1}]\to U_{l(0,j_0)}\subseteq \pi (C)$ and  $\alpha ^{\bar{j}[0],0}=F\circ a^{\bar{j}[0],0}:[t_0,t_{1}]\to U_{l(0,j_0)}\subseteq \pi (C).$ Let $\beta ^{\bar{j},0}=s_{l(0,j_0)}\circ \alpha ^{\bar{j},0},$  $\beta ^{\bar{j}[0],0}=s_{l(0,j_0)}\circ \alpha ^{\bar{j}[0],0}$ and $\rho =s_{l(0,j_0)}\circ \rho ^{\bar{j},1}.$ Let $\delta :{\maD}\to C$ be the vertical path such that $\delta _0=\beta ^{\bar{j},0}_1$ and $\delta _1=\mu _0.$

Let $\beta ^{\bar{j}}=\beta ^{\bar{j},0}\cdot \delta \cdot \mu \cdot \beta '$ and $\beta ^{\bar{j}[0]}=\beta ^{\bar{j}[0],0}\cdot \rho ^{-1}\cdot \delta \cdot \mu \cdot \beta '.$ Then $\beta ^{\bar{j}}$ and $\beta ^{\bar{j}[0]}$ are lifting of $\alpha ^{\bar{j}}$ and $\alpha ^{\bar{j}[k]}$ respectively after being modified  by constants.

Since $a ^{\bar{j}[0],0}\cdot (r ^{\bar{j},1})^{-1}\sim a^{\bar{j},0}$ by Remark \ref{nrmk as}, applying $F$ we get $\alpha ^{\bar{j}[0],0}\cdot (\rho ^{\bar{j},1})^{-1}\sim \alpha ^{\bar{j},0}$ and applying $s_{l(0,j_0)}$ we get $\beta ^{\bar{j}[0],0}\cdot \rho ^{-1}\sim \beta ^{\bar{j},0}.$ Therefore, $\beta ^{\bar{j}}\sim \beta ^{\bar{j}[0]}.$

For all the other cases the argument is similar.\\

Now if $\maI =\maI _1\wedge \ldots \wedge \maI _k$ with each $\maI_i$ a basic $d$-interval apply the previous process to each $F_{|\maI _i\times [c,d]},$ $\gamma _{|\maI _i}$ and $\sigma _{|\maI _i}$ to get $\gamma '^i,$ $\sigma '^i,$ $\beta ^i$ and $\tau ^i$ such that $\beta ^i\sim \tau ^i.$ If needed use (1) to replace the $\beta ^i$'s so that they patch together and similarly for the $\tau ^i$'s. Patch all these to get $\gamma ',$ $\sigma ',$ $\beta $ and  $\tau $ and the result follows from transitivity and compatibility with concatenation of $\sim $  (Fact \ref{fact hom and conc}).

Now if $\maJ =\maJ _1\wedge \ldots \wedge \maJ _k$ with each $\maJ_j$ a basic $d$-interval apply the previous process to each $F_{|\maI \times \maJ_j}$ and conclude by the transitivity of $\sim .$ 
\qed \\

The main consequence of Lemma \ref{lem lifting paths} is the following: \\

\begin{lem}\label{lem pi of Jcell}
 Let $C\subseteq M^{n}$ be a cell. Then:
\begin{enumerate}
\item
 $C$ is definably path connected. In fact there is a uniformly definable family of definable paths connecting a given fixed point in $C$ to any other point in $C.$
 \item
 $C$ is definably simply connected, i.e $\pi _1(C)=1.$
 \end{enumerate}
\end{lem}

\pf
(1) The proof is by complete induction on the dimension $n$ of the ambient space. If $n=0$, the space is reduced to a single point and the result follows trivially. Assume that the result holds for every $k<n.$ In order to show that it still holds for $n$ we proceed by a new induction on the definition of cells. The zero-dimensional
case is immediate. If $C=\Gamma(f) \subseteq M^{n+1}$ is the graph of a definable continuous function $f:B \to M$, where $B \subseteq M^n$ is a cell, then the projection of $C$ onto $B$ is a definable homeomorphism  and the result follows by the induction hypothesis. It remains to be considered the following case
\[
C = \left\{(x,y)\in B\times M : f(x)<y<g(x) \right\}
\]
where $f, g:B\to M$ are definable continuous maps, $B\subseteq M^n$ is a cell  and $f<g.$

If $C\subseteq M^{n+1}$ is  not open, then $\dim C<n+1$. By Remark \ref{nrmk cells perm},  $C$ is definably homeomorphic to an open cell $D \subseteq M^{\dim C}$ and the result follows by the main induction hypothesis. 

Assume that $C$ is an open cell.  By Theorem \ref{thm continuous-choice}, there is a finite cover $\{U_i:i=1,\ldots, m\}$ of $B$ by open definable subsets such that for each $i$ there is a continuous definable section $s_i:U_i\to C$  of the projection $\pi :M^{n+1}\to M^n$ onto the first $n$ coordinates (i.e.  $\pi\circ s_i={\rm id}_{U_i}$).
By Fact \ref{thm open cells}, after replacing the $U_i$'s if needed, we may assume that each $U_i$ is definably homeomorphic to an open cell in $M^n.$ So by the induction hypothesis, each $U_i$  is  definably path connected, in fact there is a uniformly definable family of definable paths connecting a given fixed point in $U_i$ to any other point in $U_i.$

For each $i$ fix $u_i\in U_i.$ Since $B$ is definably path connected, for each $i,j,$ let $\gamma ^{i,j}$ be a definable path in $B$ such that $\gamma ^{i,j}_0=u_i$ and $\gamma ^{i,j}_1=u_j.$ If needed modify each $\gamma ^{i,j}$ by constants and take by Lemma \ref{lem lifting paths} (1) a lifting $\beta ^{i,j}$ such that $\beta ^{i,j}_0=s_i(u_i)$ and $\beta ^{i,j}_1=s_j(u_j).$ Then there is a uniformly definable family of definable paths in $C$ connecting any point of $s_{i}(U_i)$ to any point of $s_j(U_j).$ Since vertical paths are uniformly definable, the same holds for $(\pi _{|C})^{-1}(U_i)$ and $(\pi _{|C})^{-1}(U_j).$ Since $C=(\pi _{|C})^{-1}(U_1)\cup \ldots \cup (\pi _{|C})^{-1}(U_m)$ the  result follows.\\

(2) The proof is again by complete induction on the dimension $n$ of the ambient space. If $n=0$, the
space is reduced to a single point and the result follows trivially. Assume that the result holds for every $k<n$. In order to show that it still holds for $n$ we proceed by a new induction on the definition of cells. The only case to be examined is the case
\[
C = \left\{(x,y)\in B\times M : f(x)<y<g(x) \right\},
\]
where $f, g:B\to M$ are definable continuous maps, $B\subseteq M^n$ is a cell  and $f<g.$  (The zero-dimensional case is immediate; if $C=\Gamma(f) \subseteq M^{n+1}$ is the graph of a definable continuous function $f:B \to M$, where $B \subseteq M^n$ is a cell, then the projection of $C$ onto $B$ is a definable homeomorphism  and the result follows by the induction hypothesis.) 

 If $C\subseteq M^{n+1}$ is not an open set, then $\dim C<n+1$. By Remark \ref{nrmk cells perm}, we have that $C$ is definably homeomorphic to an open cell $D \subseteq M^{\dim C}$ and the result follows by the main induction hypothesis. 
 
 Assume that $C$ is an open set. Let $\alpha \colon {\maI}\to C$ be a definable loop at $p \in C$. We want to show that $\alpha \approx \mathbf{c}_{\maI}^{p}.$ By the induction hypothesis, $\pi \circ \alpha \approx \mathbf{c}_{\maJ}^{\pi (p)}.$ After modifying by constants, by Lemmas \ref{lem mod by const} and \ref{lem lifting paths} (2) there are liftings $\beta $ and $\tau $ of $\pi \circ \alpha $ and $\mathbf{c}_{\maJ}^{\pi (p)}$ such that $\beta \sim \tau.$ Since after modifying by constants, $\alpha $ and $\mathbf{c}_{\maI}^{p}$ are also liftings of $\pi \circ \alpha $ and $\mathbf{c}_{\maJ}^{\pi (p)},$ by Lemmas \ref{lem unique hom lifting} and \ref{lem mod by const} and the transitivity of $\sim $ we obtain $\alpha \approx \mathbf{c}_{\maI}^{p}.$ 
\qed \\

By Fact \ref{thm open cells} and  Lemma \ref{lem pi of Jcell} we have the following which shows (P1) (a) and (P2). Compare with the corresponding results \cite[Lemma 2.9 and Proposition 3.1]{eep} in o-minimal expansions of ordered groups.\\

\begin{cor}\label{cor pi cover  def jman}
Let $X$ be a definable  manifold of dimension $n.$  Then the following hold:
\begin{enumerate}
\item
$X$ is definably connected  if and only if $X$ is definably path connected. In fact, for any definably connected definable subset $D$ of $X$ there is a uniformly definable family of definable paths in $D$ connecting a given fixed point in $D$ to any other point in $D$.
\item
$X$ has  an  admissible cover  $\{O_s\}_{s\in S}$ by  open definably  connected definable subsets such that:
\begin{itemize}
\item
$\{O_s\}_{s\in S}$ refines the definable charts of $X$;

\item
for each $s\in S$, $O_s$  is definably homeomorphic to a cell of dimension $n,$ in particular,  the o-minimal fundamental group $\pi _1(O_s)$ is trivial.\\
\end{itemize}
\end{enumerate}
\end{cor}

Finally we show  (P1) (b). For analogues compare with \cite[Section 2]{eo} in o-minimal expansions of fields or with \cite[Lemma 2.13]{eep} in o-minimal expansions of ordered groups. In all three cases the proofs are the same, they only use the fact that the domains of the corresponding definable paths and definable homotopies are definably normal. \\

\begin{lem}\label{lem ldef lifting paths}
Let $X$ and $S$ be  locally definable  manifolds with definable charts. Suppose that $p_X:X\to  S$ is a locally definable covering map. Then
the following hold.
\begin{enumerate}
\item
Let $\gamma :\maI  \to S$ be a definable path in $S.$ Let $x\in X$ be such that $p_X(x)=\gamma _{0}.$ Then there exists a unique definable path $\tilde{\gamma }:\maI \to X$ in $X$ lifting $\gamma $ such that $\tilde{\gamma } _{0}=x$.

\item
Suppose that $F: \maI \times \maJ \to S$ is a definable homotopy  between the definable paths $\gamma $ and $\sigma $ in $S$. Let $\tilde{\gamma }$ be a definable path in $X$ lifting $\gamma $. Then there exists a definable path $\tilde{\sigma }$ in $X$ lifting $\sigma $ and there exists is a unique definable lifting $\tilde{F}: \maI \times \maJ  \to X$ of $F$, which is a definable homotopy between $\tilde{\gamma }$ and $\tilde{\sigma }.$
\end{enumerate}
\end{lem}

\pf 
Let  ${\mathcal U}=\{U_{\alpha }\}_{\alpha \in I}$ be an admissible cover of $S$ by open definable subsets over which $p_X:X\to  S$ is trivial. We may assume that  ${\mathcal U}=\{U_{\alpha }\}_{\alpha \in I}$ refines the  definable charts of $S$ witnessing the fact that $S$ is a locally definable manifold with definable charts.  \\

(1) First we assume that $\maI$ is a basic $d$-interval $\langle [a,b], \langle 0_{\maI}, 1_{\maI}\rangle \rangle .$ We may also assume that the definable total order  $<_{\maI}$ on the domain $[a, b]$ of $\maI$ is $<.$ If not, the argument is similar, one just has to construct the lifting from right to left instead of from left to right.

Let $L\subseteq I$ be a finite subset such that $\gamma ([a, b])\subseteq {\bigcup_{l\in L}}U_l$. Then $[a,b]\subseteq  \bigcup_{l\in L}\gamma^{-1}(U_l)$, with  the $\gamma^{-1}(U_l)$'s open in $[a,b]$. So by Remark \ref{remark interval def normal} and the shrinking lemma (\cite[Corollary 2.12]{emp}, \cite[Chapter 6, (3.6)]{vdd}), for each $l\in L$ there is $W_l\subset [a,b]$, open in $[a,b]$ such that $W_l\subset\overline{W_l}\subset\gamma^{-1}(U_l)$ and $[a,b]\subseteq  \bigcup_{l\in L} W_l$. Therefore,  there are $a=t_0<t_1<\cdots <t_r=b$ such that for each $i=0,\dots,r-1$ we have $\gamma([t_i,t_{i+1}]) \subset U_{l(i)}$ (and $\gamma(t_{i+1})\in U_{l(i)}\cap U_{l(i+1)}$). 

Lift $\gamma_1=\gamma_{|[a,t_{1}]}$ to $\tilde{\gamma_1}=(p_{|U_{l(0)}^{i_0}})^{-1}\circ \gamma _{|[a,t_1]}$, with  $\tilde{\gamma_1}_0=x$, using the definable homeomorphism $p_{|U_{l(0)}^{i_0}}\colon U_{l(0)}^{i_0}\to U_{l(0)}$, where $U_{l(0)}^{i_0}$ is the definable connected component of $p^{-1}(U_{l(0)})$ in which $x$ lays. Repeat the process for each $\gamma_{i+1}=\gamma_{|[t_i,t_{i+1}]}$ with $\tilde{\gamma_i}(t_i)$ instead of $x$. Patch the liftings together to obtain $\tilde{\gamma }.$ 

Now if $\maI =\maI _1\wedge \ldots \wedge \maI _k$ with each $\maI_i$ a basic $d$-interval apply the previous process to lift $\gamma_1=\gamma_{|\maI _1}$ to $\tilde{\gamma_1}$, with  $\tilde{\gamma_1}_0=x$ and repeat the process for each $\gamma_{i+1}=\gamma_{|\maI _{i+1}}$ with $\tilde{\gamma_i}(1_{\maI_i})$ instead of $x$. Patch the liftings together to obtain $\tilde{\gamma }.$

Uniqueness follows (in each step) from \cite[Lemma 2.8]{eep}.\\

(2) First assume  that $\maJ$ is a basic $d$-interval $\langle [c,d], \langle 0_{\maJ}, 1_{\maJ}\rangle \rangle .$ We may also assume that the definable total order  $<_{\maJ}$ on the domain $[c, d]$ of $\maJ$ is $<$. If not, the argument is similar, one just has to construct the lifting from top to bottom  instead of from bottom  to top.

To proceed we also assume that $\maI$ is a basic $d$-interval $\langle [a,b], \langle 0_{\maI}, 1_{\maI}\rangle \rangle .$ We may furthermore  assume that the definable total order  $<_{\maI}$ on the domain $[a, b]$ of $\maI$ is $<.$ If not the argument is similar, one just has to construct the lifting from right to left instead of from left to right. 

Let $L\subseteq I$ be a finite subset such that $F([a,b]\times [c,d])\subseteq {\bigcup_{l\in L}}U_l$. Then $[a,b]\times [c,d] \subseteq \bigcup_{l\in L}F^{-1}(U_l)$, with  the $F^{-1}(U_l)$'s open in $[a, b]\times [c,d]$.  So by Remark \ref{remark interval def normal} and the shrinking lemma (\cite[Corollary 2.12]{emp}, \cite[Chapter 6, (3.6)]{vdd}), we have that  for each $l\in L$ there is $W_l\subset [a, b]\times [c,d]$, open in $[a, b]\times [a,d]$ such that $W_l\subset\overline{W_l}\subset F^{-1}(U_l)$ and $[a, b]\times[c, d] \subseteq \bigcup_{l\in L}W_l$. Now take a cell decomposition of  $[a,b]\times [c,d]$ compatible with the $W_l$'s. This cell decomposition is given by  a decomposition  $a=t_0<t_1<\cdots <t_r=b$ of $[a, b]$ together with definable continuous functions $f_{i,j}:[t_i, t_{i+1}] \to [c,d]$ for $i=0,\ldots , r-1$ and $j=0,\ldots , k_i$ such that: (i) $f_{i,0}<f_{i,1}< \ldots <f_{i, k_i}$ for $i=0,\ldots , r-1;$ (ii) $\Gamma (f_{i,0})=[t_i, t_{i+1}]\times \{c\}$ and $\Gamma (f_{i, k_i})=[t_i, t_{i+1}]\times \{d\}$ for $i=0,\ldots , r-1;$ (iii) the two-dimensional cells are of form $C_{i,j=}(f_{i,j}, f_{i, j+1})_{(t_i, t_{i+1})}.$ For each two-dimensional cell $C_{i,j}$ and each $l(i,j)$ such that $C_{i,j}\subset W_{l(i,j)}$, we have $F(\overline{C_{i,j}})\subset U_{l(i,j)}$ and  for any two-dimensional cells $C_{i,j}$ and $C_{i',j'}$ in $[a, b]\times [c,d]$,   and for each $l(i,j), l(i',j'),$ such that $C_{i,j}\subset W_{l(i,j)}$ and $C_{i',j'}\subset W_{l(i',j')}$ we also have  $F(\overline{C_{i,j}}\cap\overline{C_{i',j'}})\subset U_{l(i,j)}\cap U_{l(i',j')}$. 

Lift $F_{0,1}=F_{|\bar{C_{0,1}}}$ to $\tilde{F_{0,1}}=(p_{|U_{l(0,1)}^{i_{0,1}}})^{-1}\circ F_{|\bar{C_{0,1}}}$,  using the definable homeomorphism $p_{|U_{l(0,1)}^{i_{0,1}}}\colon U_{l(0,1)}^{i_{0,1}}\to U_{l(0,1)}$, where $U_{l(0,1)}^{i_{0,1}}$ is the definable connected component of $p^{-1}(U_{l(0,1)})$ in which $\tilde{\gamma }([t_0, t_1])$ lays. Repeat the process for each $F_{0, j+1}=F_{|\bar{C_{0, j+1}}}$ with 
$\tilde{F_{0, j}}(\Gamma (f_{0,j}))$ instead of $\tilde{\gamma }([t_0, t_1])$. Patch the liftings together to obtain $\tilde{F_{0}}:[t_0,t_1]\times [c,d]\to X$ a definable lifting of $F_{|[t_0,t_1]\times [c,d]}$ which is a definable homotopy between $\tilde{\gamma }_{|[t_0,t_1]}$ and $\tilde{\sigma }_{|[t_0,t_1]}.$ Repeat the above process again but now for each $i=1, \ldots , r-1$, starting in each case with $\tilde{\gamma }([t_i, t_{i+1}])$ and obtain the liftings $\tilde{F_{i}}:[t_i,t_{i+1}]\times [c,d]\to X$ a definable lifting of $F_{|[t_i,t_{i+1}]\times [c,d]}$ which is a definable homotopy between $\tilde{\gamma }_{|[t_i,t_{i+1}]}$ and $\tilde{\sigma }_{|[t_i,t_{i+1}]}.$ These liftings  patch together to give a  definable lifting $\tilde{F}:[a,b]\times [c,d]\to X$  of $F$ which is a definable homotopy between $\tilde{\gamma }$ and $\tilde{\sigma }.$

Now if $\maI =\maI _1\wedge \ldots \wedge \maI _k$ with each $\maI_i$ a basic $d$-interval apply the previous process to lift $F_1=F_{|\maI _1\times [c,d]}$ to $\tilde{F_1}$, with  $\tilde{F_1}(\maI _1, c)=\tilde{\gamma }(\maI _1)$ and repeat the process for each $F_{i+1}=F_{|\maI _{i+1}\times [c,d]}$ with $\tilde{\gamma }({\maI_{i+1}})$ instead of $\tilde{\gamma }(\maI _1)$. Then patch these liftings together to obtain a  definable lifting $\tilde{F}:\maI \times \maJ \to X$ of $F$ which is a definable homotopy between $\tilde{\gamma }$ and $\tilde{\sigma }.$

Now if $\maJ =\maJ _1\wedge \ldots \wedge \maJ _k$ with each $\maJ_j$ a basic $d$-interval apply the previous process to lift $F_1=F_{|\maI \times \maJ_1}$ to $\tilde{F_1}$, with  $\tilde{F_1}(\maI , 0_{\maJ_1})=\tilde{\gamma }(\maI )$ and repeat the process for each $F_{j+1}=F_{|\maI \times \maJ _{j+1}}$ with $\tilde{F_j}(\maI, 1_{\maJ _j})$ instead of $\tilde{F_1}(\maI , 0_{\maJ_1})$. To finish patch these liftings together to obtain a definable  lifting $\tilde{F}:\maI \times [c,d]\to X$ of $F$ which is a definable homotopy between $\tilde{\gamma }$ and $\tilde{\sigma }.$

As above, uniqueness follows from  \cite[Lemma 2.8]{eep}.
\qed \\



\end{subsection}

\begin{subsection}{Universal covering maps and fundamental groups}\label{subsec inv main}
As explained in Subsection \ref{subsec main prop}, from the main properties of definable paths and definable homotopies (Properties \ref{propt main new paths}) we obtain, in arbitrary o-minimal structures with definable Skolem functions, in exactly the same way as  in \cite{eep} for o-minimal expansions of ordered groups, all of the results stated below. \\

\begin{thm}\label{thm main thm}
Let $X$ be a definably connected locally definable manifold. Then:
\begin{enumerate}
\item
 there exists a universal locally definable covering map $u:U\to  X.$ Moreover, if $X$ is Lindel\"of (resp.  paracompact), then $U$ is also Lindel\"of (resp. paracompact).
 \item
 If $X$ is  Lindel\"of, then the  o-minimal fundamental group $\pi _1(X)$ of $X$  is countable. In fact, if $X$ is definable, then $\pi _1(X)$ is finitely generated.
 \end{enumerate}
\end{thm}

For similar previously known results  in special cases see  \cite{bo1}, \cite{bao2},  \cite{dk4}, \cite{EdPa},  \cite{el} and \cite{es}.

\begin{thm}\label{thm  ucovers map  ext}
Let  ${\mathbb J}$ be an elementary extension of ${\mathbb M}$ or an o-minimal expansion of ${\mathbb M}$. Let  $X$  be a definably connected locally definable manifold. Then the following hold:
\begin{enumerate}

\item
A  universal locally ${\mathbb J}$-definable covering map of $X$ is ${\mathbb J}$-definably homeomorphic to a  universal locally definable covering map of $X.$

\item
The o-minimal  fundamental group of $X$ in ${\mathbb J}$ is isomorphic to the o-minimal fundamental group of $X$ in ${\mathbb M}.$

\end{enumerate}
\end{thm}

Similarly,  we have:

\begin{thm}\label{thm  ucovers map top}
Suppose that ${\mathbb M}$ is an o-minimal  expansion of the ordered set of real numbers. Let  $X$  be a definably connected locally definable manifold. Then the following hold:
\begin{enumerate}

\item
A topological universal covering map of $X$ is topologically homeomorphic to the o-minimal universal locally definable covering map of $X.$

\item
The topological fundamental group of $X$ is isomorphic to the o-minimal fundamental group of $X.$

\end{enumerate}
\end{thm}

For previously known analogues of these invariance results in special cases see  \cite{bo1}, \cite{bao1},  \cite{bao2}), \cite{ejp2}, \cite{dk4} and \cite{dk5}. 

\begin{nrmk}
{\em
By Theorem \ref{thm  ucovers map top} when ${\mathbb M}$ is an o-minimal expansion of the ordered set of real numbers, for definably connected locally definable manifolds, the theory developed in this paper coincides with the classical theory of topological covering maps (\cite{f}). However, one should point out that, in an arbitrary  o-minimal structures ${\mathbb M},$ the theory of topological covering maps is in some sense useless. In fact in that situation, if ${\mathbb M}$ is non-archimedean, then all  definably connected locally definable manifolds are, with their natural topology, totally disconnected spaces and so have no non-trivial covering spaces. Our Theorem \ref{thm main thm} shows that it is possible to find a suitable replacement of the theory of topological covering maps which in the archimedean case coincides with the classical theory and moreover it is preserved under elementary  extensions (Theorem \ref {thm  ucovers map  ext}). \\ 
}
\end{nrmk}

In the paper \cite{Edal17} it was convenient to introduced the o-minimal $\bJ$-fundamental group which is 
the relativization of the general o-minimal fundamental group to a product of definable group-intervals. Our next goal is to show that these two kinds of o-minimal fundamental groups are isomorphic.\\

First we recall a couple of definitions. See \cite[Definition 3.1, 3.7, 3.18 and 3.19]{Edal17} (see also \cite[Definition 3.1]{epr}). 

\begin{defn}\label{defn gp-int}
{\em
A {\it definable group-interval}  $J=\langle (-b, b), 0, +,  <\rangle $ is an open interval $(-b, b)\subseteq M$, with $-b<b$ in $M\cup \{-\infty, +\infty \},$ together with a binary partial continuous definable operation $+:J^2\to J$ and an element $0\in J$, such that:
\begin{itemize}
\item[(i)]
 $x+y=y+x$ when defined; 
 $(x+y)+z=x+(y+z)$ when defined; 
 if  $x<y$ and $x+z$ and $y+z$ are defined then $x+z<y+z;$
\item[(ii)]
for every $x\in J,$ if $x>0,$ then the set $\{y\in J:\,\, x+y \,\, \textrm{is defined}\}$ is an interval of the form $(-b, r(x));$
\item[(iii)]
for every $x\in J,$ we have  $\lim _{z\to 0}(z+x)=x$ and  if $x>0$ we have also $\lim _{z\to r(x)^-}(x+z)=b;$
\item[(iv)]
for every $x\in J$ there exists $-x\in J$ such that $x+(-x)=0.$
\end{itemize}

}
\end{defn}

For the rest of this subsection let $\bJ=\Pi _{i=1}^mJ_i$  be a fixed cartesian product of   definable group-intervals $J_i=\langle (-_ib_i, b_i), $ $0_i, +_i, -_i,  <\rangle $.  

\begin{defn}
{\em
$\,$
\begin{itemize}
\item
We say that a set $X$ is a {\it $\bJ$-set} if $X\subseteq \bJ;$ in particular, a {\it $\bJ$-cell} is a cell which is  a $\bJ$-set.
\item
We say that $X$ is a (locally) definable manifold with  definable {\it $\bJ$-charts} if $X$ has definable charts  $\{(U_l, \phi _l)\}_{l\leq \kappa }$ with each $\phi _l(U_l)$ a definable  $\bJ$-set.
\item
We say that a set $X$ is a {\it $\bJ$-bounded set} if  $X\subseteq \Pi _{i=1}^m[-_ic_i,c_i]$ for some $c_i>0_i$ in $J_i;$ in particular, a {\it $\bJ$-bounded cell} is a cell which is  a $\bJ$-bounded set.
\item
We say that $X$ is a (locally) definable manifold with  definable {\it $\bJ$-bounded charts} if $X$ has definable charts  $\{(U_l, \phi _l)\}_{l\leq \kappa}$ with each $\phi _l(U_l)$ a definable  $\bJ$-bounded set.\\
\end{itemize}
}
\end{defn}

\begin{defn}
{\em
$\,\,$
\begin{itemize}
\item
A basic $d$-$\bJ$-interval is a basic $d$-interval $\maI=\langle [a,b], \langle 0_{\maI}, 1_{\maI}\rangle \rangle$ with $[a, b]\subseteq J_l$ for some $l\in \{1, \ldots ,m\};$ a $d$-$\bJ$-interval is a $d$-interval $\maI=\maI _1\wedge\dotsb\wedge \maI _n$ with each $\maI _i$ a basic  $d$-$\bJ$-interval. Note that the $\maI_i$'s can be in different $J_l$'s.
\item
If $X$ is a locally definable manifold with definable $\bJ$-charts, then a definable $\bJ$-path (resp.  constant definable $\bJ$-path, or definable $\bJ$-loop) is a definable path (resp. constant definable path or definable loop) $\alpha :\maI \to  X$ with $\maI$ a $d$-$\bJ$-interval;  $X$ is definably $\bJ$-path connected if for every $u,v$ in $X$ there is a definable $\bJ$-path $\alpha :\maI \to   X$ such that $\alpha _0=u$ and $\alpha _1=v.$ 
\item
If $X$ and $Y$ are locally definable manifolds with definable $\bJ$-charts, then two definable continuous maps $f,g:Y\to X$ are definably $\bJ$-homotopic, denoted $f\sim _{\bJ} g,$ if there is a definable homotopy $F(t,s):Y\times \maJ\to   X$  between $f$ and $g$ with $\maJ$ a $d$-$\bJ$-interval;  two definable $\bJ$-paths $\gamma :\maI \to   X$, $\delta :\maJ\to   X$, with $\gamma _0=\delta _0 $ and $\gamma _1=\delta _1$, are definably $\bJ$-homotopic, denoted $\gamma \approx _{\bJ}\delta ,$ if there  are $d$-$\bJ$-intervals $\maI'$ and~$\maJ '$
such that $\maJ '\wedge\maI=\maJ \wedge\maI'$, and there is a definable $\bJ$-homotopy 
$${\mathbf c}_{\maJ '}^{\gamma_0}\cdot \gamma
\sim _{\bJ}
\delta\cdot  {\mathbf c}_{\maI'}^{\delta_1}$$
fixing the end points. 
\end{itemize}
}
\end{defn}

The results mentioned in Subsection \ref{subsection new fund grp} for the relations $\sim $ and $\approx $ hold also for $\sim _{\bJ}$ and $\approx _{\bJ}$ respectively. 

\begin{defn}
{\em
Let $X$ be a locally definable manifolds with definable $\bJ$-charts,  $e_X\in X$ and $x_0, x_1\in X$. Let ${\mathbb P}^{\bJ}(X, x_0, x_1)$ denote the set of all  definable $\bJ$-paths in $X$ that start at $x_0$ and end at $x_1$ and let  ${\mathbb L}^{\bJ}(X, e_X)$ denotes the set of all definable $\bJ$-loops that start and end at a fixed  element $e_X$ of $X$ (i.e. ${\mathbb L}^{\bJ}(X,e_X)={\mathbb P}^{\bJ}(X,e_X,e_X)$). Then the restriction of  $\approx _{\bJ}$   to ${\mathbb P}^{\bJ}(X, x_0, x_1)\times {\mathbb P}^{\bJ}(X, x_0, x_1)$ is an equivalence relation on ${\mathbb P}^{\bJ}(X, x_0, x_1)$ and
$$\pi _1^{\bJ}(X, e_X):=
\raisebox{1ex}{$\mathsurround=0pt\displaystyle {\mathbb L}^{\bJ}(X,e _X)$}
\Big/ \raisebox{-1ex}{$\mathsurround=0pt\displaystyle \approx _{\bJ}$}
$$
 is a group, the {\it o-minimal $\bJ$-fundamental group} $\pi _1^{\bJ}(X, e_X)$ of $X,$ with group operation given by $[\gamma ][\delta ]=[\gamma \cdot \delta ]$ and identity the class a of constant $\bJ$-loop at $e_X.$ Moreover, if $f:X\to  Y$ is a locally definable continuous map between two  locally definable manifolds with definable $\bJ$-charts with $e_X\in X$ and $e_Y\in Y$ such that $f(e_X)=e_Y$, then we have an induced homomorphism $f_*:\pi _1^{\bJ}(X, e_X)\to  \pi _1^{\bJ}(Y, e_Y):[\sigma ]\mapsto [f\circ \sigma ]$ with the usual functorial properties. \\
}
\end{defn}

As usual for a definably $\bJ$-path connected locally definable manifold $X$ with definable $\bJ$-charts if there is no need to mention a base point $e_X\in X$, then by Fact \ref{fact pi1 and x and connected} (1), we may denote $\pi _1^{\bJ}(X,e_X)$ by $\pi _1^{\bJ}(X)$.\\

For ${\bJ}$-definable paths and the o-minimal ${\bJ}$-fundamental group we also have the corresponding properties (P1) and (P2).  See \cite[Corollary 3.21]{Edal17} for (P1) (a) and (P2) and see \cite[Lemma 3.23]{Edal17} for (P1) (b). Thus, just like in Theorems \ref{thm  ucovers map  ext} and \ref{thm  ucovers map  top}, we can use these properties in the two setting to prove the following, which was conjectured in the paper \cite{Edal17}: \\

\begin{thm}\label{thm  ucovers map  gpint}
Let  $\bJ=\Pi _{i=1}^mJ_i$  be a cartesian product of   definable group-intervals. Let  $X$  be a definably connected locally definable manifold with definable $\bJ$-charts. Then the following hold:
\begin{enumerate}

\item
 there exists a universal locally definable covering map $w:W\to  X$ where $W$ is a locally definable manifold with definable $\bJ$-charts. Moreover this locally definable covering is definably homeomorphic to a  universal locally definable covering map of $X.$

\item
The o-minimal  $\bJ$-fundamental group of $X$  is isomorphic to the o-minimal fundamental group of $X$ in ${\mathbb M}.$

\end{enumerate}
\end{thm}

\pf
By properties (P1) and (P2) in $\bJ$ and the proof of \cite[Theorem 1.2]{eep}, $X$ has  a  universal locally definable covering map  $w:W\to X$ where $W$ is a definably connected,  locally definable manifold with definable $\bJ$-charts,   $e_W\in W$ and $w(e_W)=e_X$. In particular, by \cite[Remark 3.8]{eep}, we  have $\pi _1^{\bJ}(W, e_W)=1.$

By Theorem \ref{thm main thm}, let $u:U\to X$ be a universal locally definable map with $U$ definably connected, $e_U\in U$ and $u(e_U)=e_X$. By 
\cite[Remark 3.8]{eep},  $\pi _1(U, e_U)=1.$ Also  there exists  a locally definable covering map $q:U\to  W$ such that
\[\xymatrix{
  U  \ar[dr]_{u}\,\,\, \ar[r]^{q } & \,\,\,\,\,\, W \ar[d]^{w}  \\
   & X
}\]
is a commutative diagram and $q(e_{U})=e_{W}$. 

Since $X$ has definable $\bJ$-charts, the same holds for $U$ (we can refine the charts of $U$ using  the admissible cover given by (P2) in $\bJ$). Therefore, $u:U\to  X$ and $q:U\to  W$ are also  locally definable covering maps in $\bJ.$ Since $\pi _1^{\bJ}(U, e_U)\simeq q_*(\pi _1^{\bJ}(U, e_U))\leq \pi _1^{\bJ}(W, e_W)=1$ (\cite[Corollary 2.17]{eep}), by \cite[Remark  3.8]{eep}, $u:U\to X$ is a universal locally definable covering map in $\bJ.$ Therefore, $u:U\to X$ and $w:W\to  X$ are locally definably homeomorphic (actually in $\bJ$) as required.

Now note that the group $\Aut (U/X)$ of locally definable homeomorphisms $\phi :U\to U$ such that $u=u\circ \phi$, is the same in ${\mathbb M}$ and in $\bJ.$ By \cite[Theorem 3.9]{eep} in ${\mathbb M}$ and in $\bJ$ respectively, we have $\pi _1(X, e_X)\simeq \Aut (U/X)$ and $\pi _1^{\bJ}(X, e_X)\simeq \Aut (U/X)$. Therefore, $\pi _1(X, e_X)\simeq \pi _1^{\bJ}(X, e_X).$
\qed

\end{subsection}

\begin{subsection}{The monodromy}\label{subsection monodromy}
Recall that if $X$ is  a locally definable manifold, then $X$ is equipped with the o-minimal site $X_{\de}$ given by: (i) the category $\op (X_{\de})$ of open definable subsets of $X$ with morphisms being inclusions; (ii) the Grothendieck topology such that for $U\in \op (X_{\de})$,  a collection $\{U_j\}_{j\in J}$ of objects of $\op (X_{\de})$ is an admissible cover of $U$ if it admits  a finite subcover. 

If $\mathsf C$ is any category  admitting projective and inductive limits and satisfying the IPC property (see \cite[Definition 3.1.10]{ks2} for more details), 
then the category of $\mathsf C$-pre-sheaves on the o-minimal site $X_{\de},$ denoted  $\psh _{\mathsf C}(X_{\de}),$ is  the category $\fct (\op (X_{\de})^{{\rm op}}, \mathsf C)$ of contravariant functors
\begin{eqnarray*}
& & \F: \op(X_{\de})\to \mathsf C\\
& & \,\,\,\,\,\, \,\,\,\,\,\, \,\,\,\,\,\, U\mapsto \F (U)\\
& & \,\,\,\,\,\, \,\,\,\,\,\, \,\,\,\,\,\, (V\subset U)\mapsto (\F (U)\to \F (V))\\
& &  \,\,\,\,\,\,\,\,\,\,\,\,\,\,\,\,\, \,\,\,\,\,\,\,\,\,\,\,\,\,\,\,\,\,\,\,\, \,\,\,\,\,\,\,\,\,\, \,\,\,\,\,\, s\mapsto s_{|V}
\end{eqnarray*}
from  $\op (X_{\de})$ to $ \mathsf C$ with morphisms being natural transformations of such functors. The category of $\mathsf C$-sheaves on the o-minimal site $X_{\de},$  denote $\sh _{\mathsf C}(X_{\de}),$ is the full subcategory of $\psh _{\mathsf C} (X_{\de})$ whose objects satisfy the following gluing conditions: for every
$U\in \op(X_{\de})$ and every admissible cover $\{U_j\}_{j\in J}$ of $U$ we have:

\begin{itemize}
\item
if $s, t\in \F (U)$ and $s_{|U_j}=t_{|U_j}$ for each $j$, then $s=t$;

\item
if $s_j\in \F (U_j)$ are such that $s_j=s_k$ on $U_j\cap U_k$ then they glue to $s\in \F (U)$ (i.e. $s_{|U_j}=s_j$).
\end{itemize}

If $V\in  \op (X_{\de})$, a ${\mathsf C}$-sheaf $\F$ on $V_{\de}$ is constant if it is isomorphic to the ${\mathsf C}$-sheaf $C_V$ on $V_{\de}$ associated to the ${\mathsf C}$-pre-sheaf sending every $W\in \op (V_{\de})$ to a fixed  $C \in \ob {\mathsf C}.$ We denote by $\csh _{\mathsf C}(X_{\de})$ the category  of constant $\mathsf C$-sheaves on the o-minimal site $X_{\de}$ on $X$. We denote by $\lcsh _{\mathsf C}(X_{\de})$ the  category of locally constant $\mathsf C$-sheaves on the o-minimal site $X_{\de}$ on $X$.  By definition, this means that, $\F \in \ob \lcsh _{\mathsf C}(X_{\de})$ if there exists an admissible cover  $\{U_j\}_{j \in J}$ of $X$ by open definable subsets such that the restriction $\F _{|U_j}$ is a constant ${\mathsf C}$-sheaf on $U_{j \, \de}$ for each $j \in J.$ (For further details on the theory of o-minimal sheaves we refer to, for example, \cite{ejp1} and \cite{ep1}).\\

Just like in \cite{eep}, from Properties \ref{propt main new paths} we obtain the  monodromy equivalence for locally constant o-minimal sheaves:

\begin{thm}\label{thm mu j equiv intro}
Then the  monodromy functor
$$\mu : \lcsh _{\mathsf C} (X_{\de})\to \fct (\Pi _1(X), \mathsf C  )$$
is an equivalence between the  category  of locally constant $\mathsf C$-sheaves on the o-minimal site $X_{\de}$ on $X$ and the category of representations of the o-minimal fundamental groupoid $\Pi _1(X)$ of $X$ in $\mathsf C.$
\end{thm}

Note that when $X$ is definably connected and $x\in X,$ then $\fct (\Pi _1(X), \mathsf C  )$ is the category of representations of the o-minimal fundamental group $\pi _1(X, x)$ of $X$ in $\mathsf C.$

Taking for $\mathsf C$ the category of $\pi _1(X,x)-$sets or of $G$-torsors, we obtain from Theorem \ref{thm mu j equiv intro}  classification results for locally definable covering maps, the  o-minimal Hurewicz and Seifert - van Kampen theorems just like in the case of o-minimal expansions of ordered groups in \cite[Subsection 4.3]{eep}.  Analogues of the o-minimal Hurewicz and Seifert - van Kampen theorems for definable sets in o-minimal expansions of fields were proved before in \cite{eo} and \cite{bo1} respectively.\\

\end{subsection}
\end{section}

\end{document}